\documentclass[11pt, reqno, a4paper]{amsart}
\usepackage[margin=2.5cm]{geometry}
\usepackage[margin=0.6cm]{caption}
\usepackage[foot]{amsaddr}

\usepackage{amsmath}
\usepackage{amsfonts}
\usepackage{array}
\usepackage{graphicx}
\usepackage{epsfig}
\usepackage{color}
\usepackage{enumitem}  
\usepackage{epstopdf}
\usepackage{stmaryrd}
\usepackage{yfonts}
\usepackage{diagbox}
\usepackage[utf8]{inputenc}
\usepackage{booktabs}  
\usepackage{capt-of}
\usepackage{multirow}
\usepackage[displaymath, mathlines]{lineno}
\usepackage{standalone}
\usepackage[colorinlistoftodos]{todonotes}

\usepackage{cleveref}
\usepackage{comment}

\usepackage{floatrow}
\newfloatcommand{capbtabbox}{table}[][\FBwidth]

\def\half{1/2}
\def\ud{u_{\odot}}

\def\udhat{\hat{u}_{\odot}}
\def\qdhat{\hat{q}_{\odot}}

\def\vd{v_{\odot}}
\def\Vd{V_{\odot}}

\def\ld{p_{\odot}}
\def\Ld{Q_{\odot}}

\def\md{q_{\odot}}

\def\eps{\varepsilon}

\def\TT{\boldsymbol T}

\def\R{\mathbb{R}}

\def\trace{{\mathcal{T}_\Gamma}}
\def\mtrace{{\overline{\mathcal{T}}_\Lambda}}
\def\tracel{{\mathcal{T}_\Lambda}}

\newcommand{\jump}[1]{\ensuremath{[\![#1]\!]} }

\newcommand{\avrd}[1]{\overline{\overline{#1}}}
\newcommand{\avrc}[1]{\overline{#1}}

\renewcommand{\TT}{\mathcal{T}}
\newcommand{\vertiii}[1]{{\left\vert\kern-0.25ex\left\vert\kern-0.25ex\left\vert #1 
    \right\vert\kern-0.25ex\right\vert\kern-0.25ex\right\vert}}

\newtheorem{theorem}{Theorem}[section]

\newtheorem{lemma}{Lemma}[section]

\begin{document}

\title[Prenconditioners for 3\textit{d}-1\textit{d} coupling with Lagrange multipliers]{Robust Preconditioning of mixed-dimensional PDEs on 3\textit{d}-1\textit{d} domains coupled with Lagrange multipliers}
\author[N. Dimola]{Nunzio Dimola$^\dag$}
\address{$^\dag$MOX, Department of Mathematics\thanks{Submitted to \textit{Quantitative approaches to microcirculation: mathematical models, computational methods, measurements and data analysis (SEMA/SIMAI Collection on Microcirculation)}}, Politecnico di Milano, Italy}
\email{nunzio.dimola@polimi.it}

\author[M. Kuchta]{Miroslav Kuchta$^*$}
\address{$^*$Simula Research Laboratory, Oslo, Norway}
\email{miroslav@simula.no}

\author[K.-A. Mardal]{Kent-Andre Mardal$^\ddag$}
\address{$^\ddag$Department of Mathematics, University of Oslo, Norway}
\email{kent-and@math.uio.no}

\author[P. Zunino]{Paolo Zunino$^\dag$}
\email{paolo.zunino@polimi.it}




%
%
\maketitle


\begin{abstract}
In the context of micro-circulation, the coexistence of two distinct length scales - the vascular radius and the tissue/organ scale - with a substantial difference in magnitude, poses significant challenges. To handle slender inclusions and simplify the geometry involved, a technique called \textit{topological dimensionality reduction} is employed, which suppresses manifold dimensions associated with the smaller characteristic length. However, the resulting discretized system's algebraic structure presents a challenge in constructing efficient solution algorithms.
This chapter addresses this challenge by developing a robust preconditioner for the $3d$-$1d$ problem using the operator preconditioning technique. Robustness of the preconditioner is demonstrated with respect to problem parameters, except for the vascular radius. The vascular radius, as demonstrated, plays a fundamental role in mathematical well-posedness of the problem and the preconditioner's effectiveness.
\end{abstract}

\section{Introduction}

The human cardiovascular system displays a diverse range of scales and characteristics, encompassing major blood vessels, arterioles, and capillaries, with diameters ranging from several cm to a few $\mu$m. In particular, when examining the microcirculation, the challenge intensifies due to the substantial disparity in scale between the diameter of small vessels ($\mu$m) and the size of the corresponding system or organ they supply (dm). To cope with this disparity, intricate vascular networks that occupy space are needed, leading to complex geometric structures.

Due to the inherent complexities involved, simulating the flow throughout the entire cardiovascular system is not practical or applicable. However, it is crucial to adopt a comprehensive approach that incorporates interactions between different system components by integrating models operating at various levels of detail. This approach is commonly referred to as \textit{geometrical multi-scale} or \textit{hybrid-dimensional modeling}.

Two types of hybrid dimensional models can be considered: sequential and embedded. In an embedded multiscale model, components of varying levels of detail are integrated within the same domain. A prime example is the micro-circulation, where a complex vascular network comprising arterioles, capillaries, and small veins exists within the biological tissue.

From the modeling standpoint, such ideas have appeared in the past three decades, (at least), for modeling wells in subsurface reservoirs in \cite{Peaceman1978183,Peaceman1983531} and for modeling microcirculation in \cite{Blake1982173,Fleischman1986141,Flieschman1986145,Secomb20041519,10.1371/journal.pcbi.1008584}. 
A similar approach has been recently used to model soil/root interactions \cite{Koch2018}.
However, these application-driven seminal ideas were not followed by a systematic
theory and rigorous mathematical analysis. At the same time the models introduce additional mathematical complexity,
in particular concerning the functional setting for the solution as it involves coupling PDEs on domains with
high dimensionality gap. Such mathematical challenge has recently attracted the attention of many researchers. The sequence of works by~\cite{DAngelo,d2012finite,d2008coupling}, 
followed by \cite{laurino_m2an,Koppl2018953,kuchta2021analysis},
have remedied the well-posedness by weakening the regularity assumptions that define a solution.

Nevertheless, the mathematical understanding of embedded, mixed-dimensional problems is not sufficient to successfully apply these models to realistic problems. An essential difficulty to overcome is the development of efficient numerical solvers that can handle the large separation of spatial scales between the domains of the equations and the subsequent geometrical complexity of the vascular network.

Actually, the study of the interplay between the mathematical structure of the problem and solvers, as well as preconditioners
for its discretization is still in its infancy. The results presented in \cite{kuchta2016preconditioners} for the solution of
$1d$ differential equations embedded in $2d$, and more recently extended to the $3d$-$1d$ case in \cite{KMM2}, have paved the way,
but a lot more has to be understood. This work proceeds along this direction, with the aim to discuss the main mathematical
challenges at the basis of the development of optimal solvers of mixed-dimensional $3d$-$1d$ problems coupled by Lagrange multipliers.

\begin{figure}[t!]
\begin{center}
\includegraphics[width=0.5\textwidth]{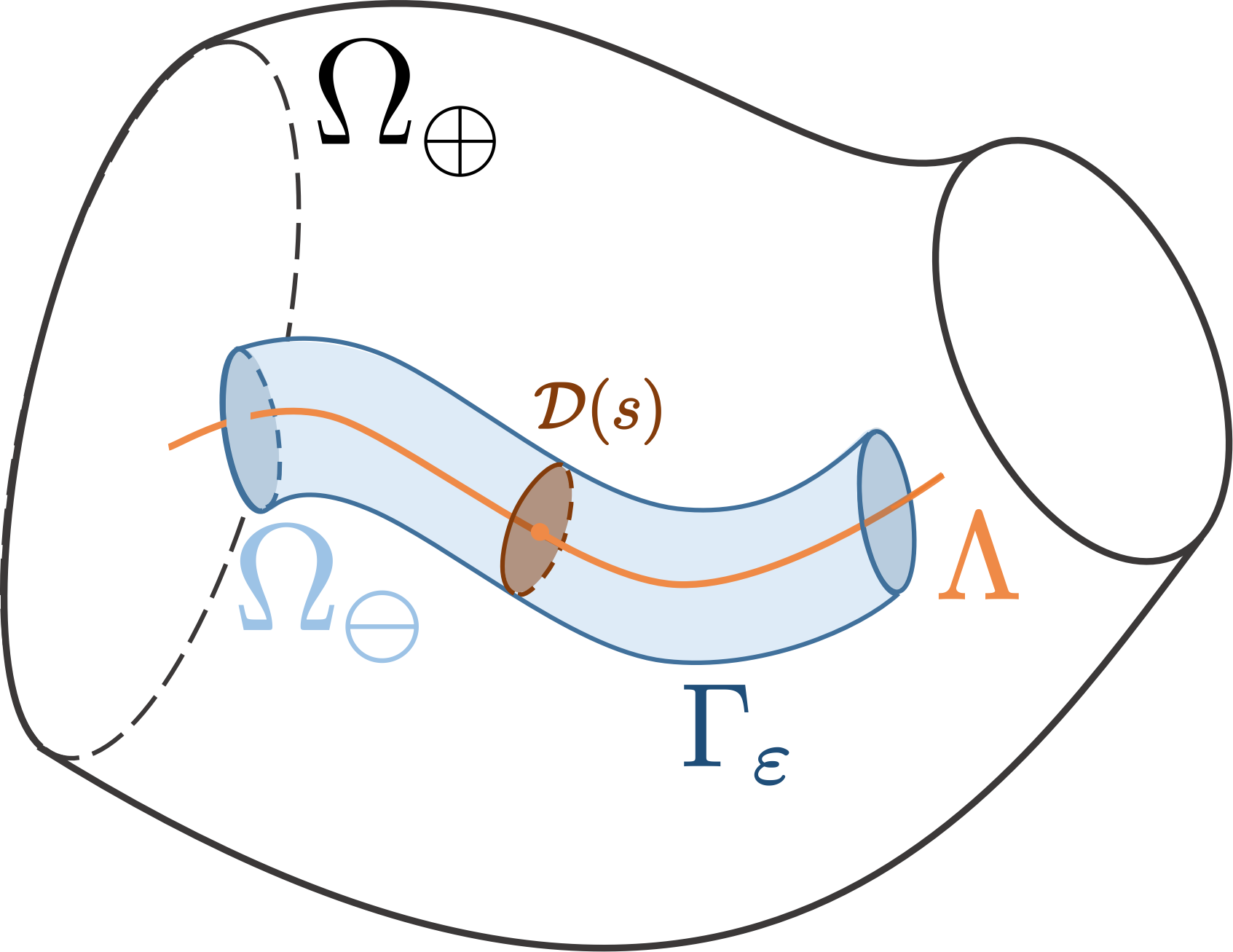}
\end{center}
\vspace{-5pt}
\caption{
  Geometrical setting considered in this paper, and the associated notation.
  Domain $\Omega_{\ominus}$ is assumed to be a slender generalized cylinder characterized
  by the centerline $\Lambda$ parameterized by the arclength coordinate $s$ and $2d$ transversal
  cross section $\mathcal{D}(s)$.
}
\label{fig1}
\end{figure}

We model the case of one single small vessel by the domain $\Omega\subset \R^3 $ that is an open,
connected and convex set that can be subdivided in two parts,  $\Omega_{\ominus}$ and
$\Omega_{\oplus}=\Omega\setminus\overline{\Omega}_{\ominus}$ representing respectively
the vessel and the tissue. Let  $\Omega_{\ominus}$ be a \emph{generalized cylinder}
that is the swept volume of a two dimensional set, $\partial\mathcal{D}$, moved along a curve,
$\Lambda$, in the three-dimensional domain, $\Omega$, see \Cref{fig1} for an illustration.
For dimensional reduction of the vessel domain we finally assume that $\lvert \mathcal{D} \rvert \ll \lvert \Lambda \rvert $.

With the purpose of developing robust preconditioners for the \emph{3d-1d} problem arising
from microcirculation, let us consider a prototype problem
that originates from coupling of two $3d$ diffusion equations:
\begin{equation}
\label{eq:pde3d}
\begin{aligned}
- \kappa_\oplus \Delta u_\oplus &= f &\mbox{ in } \Omega_\oplus,\\
- \kappa_\ominus \Delta u_\ominus &= f &\mbox{ in } \Omega_\ominus,\\
u_\oplus - u_\ominus &=  g &\mbox{ on } \Gamma_{\eps},\\
\kappa_\oplus \nabla u_\oplus \cdot {\nu}_{\Gamma_{\eps}}
- \kappa_\ominus \nabla u_\ominus \cdot {\nu}_{\Gamma_{\eps}}  &= 0 &\mbox{ on } \Gamma_{\eps},\\
u_\oplus,\ u_\ominus &= 0 &\mbox{ on } \partial\Omega. 
\end{aligned}
\end{equation}
Here, $u_\oplus$ and $u_\ominus$ are the unknowns, and $\kappa_\oplus>0$ and $\kappa_\ominus>0$ are the diffusivities in the two
different domains which we shall in the following assume to be constant for the sake of simplicity.
Further, $\Gamma_{\eps}= \partial\Omega_\oplus\cap\partial\Omega_\ominus$
where $\eps=\text{diam}\mathcal{D}$ denotes the diameter of the 2$d$ transversal
cross sections of $\Omega_{\ominus}$. Finally, $f$ is a forcing term defined in the whole $\Omega$.
We remark that the coupling condition $u_\oplus - u_\ominus =  g$ on $\Gamma_{\eps}$ is of Dirichlet
type. In this sense our starting point differs from the more commonly studied 
problem (see the seminal paper \cite{d2008coupling}) which considers the Robin coupling,
i.e. modeling the flux at the interface as
$-\kappa_{\oplus}\nabla u_\oplus\cdot \nu_{\Gamma_{\eps}} = \gamma (u_\oplus-u_\ominus)$ on $\Gamma_{\eps}$.
However, viewing the latter as a perturbation of the Dirichlet condition, our setting is relevant
also for the Robin coupling, see \cite{braess1996stability}.

By applying a suitable model reduction strategy to \eqref{eq:pde3d}, which exploits the following transverse averages as already described in \cite{laurino_m2an},
\begin{gather*}
\avrc{w}(s) = |\partial\mathcal{D}(s)|^{-1} \int_{\partial\mathcal{D}(s)} w d\gamma \quad \mbox{and} \quad 
\avrd{w}(s) = |\mathcal{D}(s)|^{-1} \int_{\mathcal{D}(s)} w d\sigma, 
\end{gather*}
we obtain the $\emph{3d-1d}$ coupled problem of the form
\begin{equation}
\label{eq:pde}
\begin{aligned}
  -\kappa\Delta u + \eps p_{\odot} \delta_{\Lambda} &= f &\mbox{ in } \Omega,\\
-\eps^2\kappa_{\odot} d_s^2 \ud  - \eps p_{\odot} &= \eps^2 \avrd{f} &\mbox{ on } \Lambda,\\
\eps(\mtrace{u} - \ud)  &=  g &\mbox{ on } \Lambda,\\
u &= 0 &\mbox{ on }\partial\Omega,\\
\ud &=0 &\mbox{ on }\partial\Lambda.\\
\end{aligned}
\end{equation}
Here $\kappa=\kappa_\oplus$ in $\Omega_\oplus$ and $\kappa_\odot = \avrd{\kappa}_\ominus$ on $\Lambda$ being the
centerline of $\Omega_\ominus$.
Throughout the paper we will use the subscript $\odot$, e.g. $\ud$, for functions on
$\Lambda$ (and, in general, on domains with topological dimension smaller than that of the ambient space)
while functions on $\Omega$ are denoted, as usual, by italic lower case letters.  
Here, the primal unknowns are $u$ and $\ud$, while the unknown ${p_{\odot}}$ is the
Lagrange multiplier used to enforce the coupling of $u$ and $\ud$. Furthermore, $\delta_\Lambda$ is
a Dirac delta function. The operator obtained from a combination of the average operator $\avrc{(\cdot)}$
with the trace on $\Gamma_{\eps}$ will be denoted with $\mtrace = \avrc{(\cdot)} \circ \trace$, as it maps
functions on $\Omega$ to functions on $\Lambda$. Note that $\mtrace$ thus implicitly depends on $\eps$.

Because of the high dimensional gap between $3d$ and $1d$, and the presence of singularities in the
solution $u$ due to the $1d$ coupling, the well-posedness of the problem \eqref{eq:pde} is 
challenging~\cite{d2008coupling, d2012finite}. In the context of Robin coupling conditions at the $3d$-$1d$ interface,
the singularity is addressed in a number of works e.g. by proposing a splitting
strategy where the singularity is captured by Green's functions 
\cite{gjerde2019singularity}, regularizing the problem by distributing the
singular source term on the coupling surface or in the bulk \cite{KVWZ, koch2020modeling} or deriving 
error estimates in tailored norms \cite{koppl2014optimal, d2012finite}.
For these formulations efficient solution algorithms have been developed utilizing
algebraic multigrid \cite{hu2023effective, cerroni2019mathematical} or fast
Fourier transform for preconditioning \cite{linninger2023mesh}.

While the Robin coupling typically leads to elliptic problems, enforcing the Dirichlet
coupling in \eqref{eq:pde} by Lagrange multiplier yields a saddle point system. 
Here, as recently shown in \cite{kuchta2021analysis}, the high $3d$-$1d$ dimensional
gap requires the multiplier to reside in fractional order Sobolev spaces on $\Lambda$ in order
to obtain well-posedness. However, scalable solvers for the resulting linear systems are
currently lacking. Here we aim to address this issue by constructing robust
block-diagonal preconditioners for \eqref{eq:pde} within the framework
of operator preconditioning \cite{mardal2011preconditioning}. In particular, to obtain
robust estimates we establish well-posedness of the coupled problem in weighted, fractional
Hilbert spaces.

The structure of the chapter is outlined as follows. We start with the analysis of the coupled problem in \Cref{sec:3d1d_analysis}. Next, we investigate the performance of the resulting preconditioner in \Cref{sec:3d1d_precond}, leading us to the observation that the parameter $\eps$ in \eqref{eq:pde} assumes a distinctive role beyond that of a standard material parameter. This unique role is further examined in  \Cref{sec:3d1d_issue}, where we establish a close connection between $\eps$ and the well-posedness of the coupled problem. Finally, we provide our conclusions in \Cref{sec:conclude}.

\section{Mathematical formulation of the \emph{3d-1d} coupled problem}\label{sec:3d1d_analysis}



Before delving into the precise mathematical formulation of problem \eqref{eq:pde}, we will introduce some fundamental notation. Let $D$ represent a bounded domain in $\mathbb{R}^d$, where $d=1,2,3$. In this context, $L^2(D)$ denotes the space of functions that are square integrable over $D$, while $H^k(D)$ refers to the space of functions possessing $k$ derivatives in $L^2$. We denote the closure of $C^\infty_0(D)$ in $H^k(D)$ as $H^k_0(D)$. When the context is clear, we may simplify the notation by using $H^k$ instead of $H^k(D)$. For a Hilbert space $X = X(D)$ with its corresponding dual space $X'$, the norm and associated inner products are denoted as $\lVert \cdot \rVert_X$ and $(\cdot, \cdot )_X$, respectively. The duality pairing is represented by $\langle \cdot, \cdot \rangle_X$, and if necessary, we will explicitly indicate the underlying domain by using $(\cdot, \cdot )_D$ and $\langle \cdot, \cdot \rangle_D$. With a slight abuse of notation, the dual of $H^k$ is $H^{-k}$ for $k\in \mathbb{R}$. Similarly, $H^{-k}_0$ is the dual of $H^k_0$.  We will use scaled and intersected inner products, norms and Hilbert spaces for our analysis of robustness. For instance $\eps H^1_0(D)$ is defined in terms of the inner product
\[
(u,v)_{\eps H^1_0(D)} = \eps^2 (\nabla u, \nabla v)_{L^2(D)} .
\]
Notice that with this definition, the norm scales linearly in $\eps$, i.e., 
\[
\|u\|_{\eps H^1_0} = \eps \|u\|_{H^1_0}.
\]
The intersection of two Hilbert spaces $X$ and $Y$ is denoted $X\cap Y$ and has an inner product
\[
(x, y)_{X\cap Y} = (x, y)_X + (x, y)_Y . 
\]

We can now formulate precisely the weak formulation of problem \eqref{eq:pde} that will be considered in this paper.  
The problem reads: Find $u\in V, \ud\in \Vd, \ld\in \Ld$ such that 
\begin{equation}
\label{eq:3D1Dvar}
\begin{aligned}
a([u, \ud], [v, \vd]) + b([v, \vd], \ld) &= \langle f, [v, \vd]\rangle_{\Omega} &\quad \forall [v, \vd] \in [V, \Vd] \\ 
b([u, \ud], \md) &= \langle g,\md\rangle_{\Lambda} &\quad \forall \md \in \Ld. 
\end{aligned}   
\end{equation}
with
\begin{align*}
a([u, \ud], [v, \vd])&= \kappa (\nabla u, \nabla )_{L^2(\Omega)} + \kappa_{\odot}\eps^2(d_s \ud, d_s \vd)_{L^2(\Lambda)},
\\
b([v, \vd], \md)&= \eps \langle \mtrace v - \vd , \md \rangle_\Lambda
\end{align*}
and $V=\kappa^{\half}H_0^1(\Omega)$, $\Vd=\kappa_{\odot}^{\half}\eps H^1_0(\Lambda)$, and  $\Ld=\kappa^{-\half}\eps H^{-\frac 1 2}(\Lambda)\cap  \kappa_{\odot}^{-\half}H^{-1}(\Lambda)$.
We remark that a main difference between \cite{kuchta2021analysis} and the functional
setting considered here is that the Lagrange multiplier space $\Ld$ is a weighted
intersection space $\kappa^{-\half}\eps H^{-\frac 1 2}(\Lambda)\cap  \kappa_{\odot}^{-\half}H^{-1}(\Lambda)$
rather than $H^{-\frac 1 2}(\Lambda)$. This choice is fundamental for the development of robust
preconditioners, as will be motivated later on.

\subsection{The stability of the continuous problem}
\label{3d-1d-stability}
We show well-posedness of the variational problem \eqref{eq:3D1Dvar} in 
space $X= V \times V_{\odot} \times Q_{\odot}$ considered with its canonical norms
by applying the abstract Brezzi theory \cite{brezzi1974existence}. To this end we let
\begin{equation*}
\begin{aligned}
  &A: (V \times V_{\odot})\rightarrow(V \times V_{\odot})', &\langle A [u, u_{\odot}], [v, v_{\odot}]\rangle_{V \times V_{\odot}} = a([u, u_{\odot}], [v, v_{\odot}]),\\
  &B: (V \times V_{\odot})\rightarrow Q_{\odot}', &\langle B [u, u_{\odot}], q_{\odot} \rangle_{Q_{\odot}} = b([u, u_{\odot}], q_{\odot}),\\
\end{aligned}
\end{equation*}
such that \eqref{eq:3D1Dvar} can be equivalently stated as: Find $x=[u, \ud, \ld] \in X$ such that 
\[
\mathcal{A}x = f \text{ in } X',  \quad \mathcal{A}=\begin{pmatrix}A & B'\\ B & 0\end{pmatrix}.
\]
The operator $\mathcal{A}:X\rightarrow X'$ is then an isomorphism (equivalently \eqref{eq:3D1Dvar} is well-posed)
provided that it satisfies the four Brezzi conditions, namely, the operators $A$ and $B$ are bounded,
$A$ is coercive on the $\text{ker}B \subset V\times V_{\odot}$ and $B$ satisfies
the inf-sup condition.


Of the four Brezzi conditions required for problem \eqref{eq:3D1Dvar}, the 
coercivity and boundedness of $A$ were established in \cite{kuchta2021analysis} in the
setting required here ($V\times V_{\odot}$ are the same and only $Q_{\odot}$ has changed)  and will therefore not be repeated. Hence, we only need to 
verify the properties of the $B$ operator.

\begin{theorem}\label{th:bnb}
$b(\cdot,\cdot): [V,\Vd] \times \Ld  \rightarrow \mathbb{R}$ satisfies the following conditions: 
\begin{align}
\label{BNB1}
&b([u, \ud],\md) \le C_B \|[u,\ud]\|_{[V, \Vd]} \|\md\|_{\Ld}, & u\in V, \ud\in \Vd, \md\in \Ld,\\[10pt] 
\label{BNB2}
&\sup_{[u, \ud]\in [V,\Vd]} \frac{b([u, \ud],\md)}{\|[u, \ud]\|_{[V, \Vd]}} \geq \beta \|\md\|_{\Ld}, & \md \in \Ld
\end{align}
with positive constants $C_B$, $\beta$. 
\end{theorem} 

\begin{proof}
  We recall that owing to Dirichlet boundary conditions for $u$, $u_{\odot}$ we
  have the equivalence between the seminorms (induced by $A$) and the full $H^1$ norms.
  In \cite{kuchta2021analysis} the following was established: 
\begin{equation*}
    \eps \|\mtrace u\|_{H^{1/2}(\Lambda)} \le C_\TT \|\TT_{\Gamma_{\eps}} u\|_{H^{1/2}(\Gamma_{\eps})} \le  C \|u\|_{H^1(\Omega)},  \quad \forall u \in H^1(\Omega),
\end{equation*}
where constants $C_\TT$ and $C$ are independent from $\eps$.
Hence, the boundedness of $b(\cdot, \cdot)$ can be obtained in a direct way as: 
\begin{align*}
b([u,& \ud], \md) = 
\eps \langle \mtrace u -   \ud, \md \rangle_{\Lambda} \\
&\le
\eps |\langle \mtrace u, \md \rangle_{\Lambda}| +  
\eps |\langle \ud, \md \rangle_{\Lambda}| \\
&\le \|\mtrace u\|_{H^{1/2}(\Lambda)} \eps\| \md\|_{H^{-1/2}(\Lambda)} +  
\eps \|\ud\|_{H^1(\Lambda)} \|\md\|_{H^{-1}(\Lambda)}  \\
&= \kappa^{1/2} \eps\|\mtrace u\|_{H^{1/2}(\Lambda)} \kappa^{-1/2}\| \md\|_{H^{-1/2}(\Lambda)} +  
\eps \kappa^{1/2}_{\odot}\|\ud\|_{H^1(\Lambda)} \kappa_{\odot}^{-1/2}\|\md\|_{H^{-1}(\Lambda)}  \\
&\le \kappa^{1/2} C\| u\|_{H^{1}(\Omega)}  \kappa^{-1/2}\| \md\|_{H^{-1/2}(\Lambda)} +  
\eps \kappa^{1/2}_{\odot}\|\ud\|_{H^1(\Lambda)} \kappa_{\odot}^{-1/2}\|\md\|_{H^{-1}(\Lambda)}  \\
&  \le C_B (\|u\|^2_V + \|\ud\|^2_{\Vd})^{1/2} \|\md\|_{\Ld}.
\end{align*}

Next, the inf-sup condition reads
\[
\sup_{[u,\ud] \in [V, \Vd]} \frac{\eps \langle \mtrace u -   \ud, \md \rangle_{\Lambda}}
    {(\|u\|^2_{\sqrt\kappa H^1(\Omega)} + \eps^2\|\ud\|^2_{\sqrt\kappa_{\odot}H^1(\Lambda)})^{1/2}}
    \ge \beta \|\md\|_{Q_{\odot}}. 
\]
We first note that $(\alpha X)'=\alpha^{-1}X'$ holds for dual of weighted
Hilbert space. In turn, let $R_1$ be the Riesz map from $\kappa_{\odot}^{-1/2} H^{-1}$ to $\kappa_{\odot}^{1/2}H^1$ 
and $R_{1/2}$ be the Riesz map from $\kappa^{-1/2} H^{-1/2}$ to $\kappa^{1/2}H^{1/2}$. 
Then $\udhat = -\frac{1}{\eps} R_1 \md$ implies that $\eps \| \udhat \|_{\sqrt{\kappa_{\odot}} H^{1}} = \|\md\|_{\kappa_{\odot}^{-1/2}H^{-1}}$.  
From \cite{kuchta2021analysis}, we also have
that there exists a harmonic extension operator $H$ such 
that for $\hat{u} = H \qdhat, \qdhat \in H^{1/2}$. Letting $\qdhat = R_{1/2} \md$ we obtain that  $\qdhat = \mtrace \hat{u}$
where $\|\qdhat\|_{\sqrt{\kappa} H^{1/2}} = \|\md\|_{\kappa^{-1/2}H^{-1/2}}$ and
\[\|\hat{u}\|_{H^1(\Omega)} \le C_{IT} \eps \|\qdhat\|_{H^{\frac12}(\Gamma_{\eps})} = \eps C_{IT} \|\md\|_{H^{-\frac12}(\Gamma_{\eps})}.\]
We note that the latter inequality, generally named the \emph{inverse trace inequality},
involves the generic constant $C_{IT}$ that may possibly depend on the domain
$\Gamma_{\eps}$ and precisely on its cross section quantified by the (small) parameter $\eps$.
Hence, we obtain 
\begin{align*}
&    \sup_{[u, \ud] \in [V, \Vd]} \frac{\eps \langle \mtrace u -   \ud, \md \rangle_{\Lambda}}
{(\|u\|^2_{\sqrt\kappa H^1(\Omega)} + \eps^2\|\ud\|^2_{\sqrt\kappa_{\odot} H^1(\Lambda)})^{1/2}} \ge 
 \frac{\eps \langle \mtrace \hat{u} -   \udhat, \md \rangle_{\Lambda}}
{ (\|\hat{u}\|^2_{\sqrt\kappa H^1(\Omega)} + \eps^2\|\udhat\|^2_{\sqrt\kappa_{\odot} H^1(\Lambda)})^{1/2}  } \\
& \ge
 \frac{\eps (\langle R_{1/2} \md, \md \rangle +    \langle \frac{1}{\eps } R_1 \md, \md \rangle  )}
{ (C_{IT} \eps^2 \|\md\|^2_{\kappa^{-1/2}H^{-1/2}(\Lambda)}+ \|\md\|^2_{\kappa_{\odot}^{-1/2} H^{-1}(\Lambda)} )^{1/2} }
\ge \beta \|\md\|_{Q_{\odot}}.
\end{align*}
\end{proof}

It is noteworthy that a proper choice of weighted-intersected Sobolev spaces has led to removal of any \emph{explicit}
dependence of the inf-sup constant $\beta$ on $\kappa$, $\kappa_{\odot}$ and $\eps$.
On the other hand, caution is needed as, through possible dependence of $C_{IT}$ on $\eps$,
the radius may influence $\beta$, and in turn robustness of the estimates in $\eps$.


\subsection{Numerical evidence about preconditioning the mixed-dimensional problems}
The exploitation of iterative solvers, such as Krylov methods, is crucial to
obtain fast solution algorithms for the discretized problems, yet, their
performance essentially depends on existence of (a practical) preconditioner
which improves spectral properties of the linear systems \cite{saad2003iterative}.
In constructing preconditioners for the $\emph{3d-1d}$-problems, a main
objective is to establish a parameter robustness preconditioner, i.e. a preconditioner 
with  performance that does not deteriote  with respect to  
discretization parameters (mesh and element sizes),  
 variations in material parameters ($\kappa>0$, $\kappa_{\odot}>0$)  as well as
the geometrical parameter $\eps$, crucial for \emph{3d-1d} problems.  As an example, in the context of \emph{3d-1d} microcirculation
model \eqref{eq:pde}, the vessel radius appears explicitly and implicitly in the \emph{3d-1d}
problem formulation, as the residue of the bridging of the three-dimensional topology to the one-dimensional one. Then, it must be taken into account.

Here we wish to construct robust preconditioners for the coupled $3d$-$1d$
problem \eqref{eq:3D1Dvar} which induces an operator equation: Find $x=[u, u_{\odot}, p_{\odot}]\in X$
such that
\begin{equation}\label{eq:coupled_3d_1d_operator}
  \mathcal{A}x = f\mbox{ in }X',\quad
\mathcal{A}
 := \begin{pmatrix}
  -\kappa\Delta & & \eps\mtrace'\\
  & -\kappa_{\odot}\eps^2\Delta_{\Lambda} & -\eps I \\
  \eps\mtrace & -\eps I & 0\\
\end{pmatrix},
\end{equation}
where, for a one-dimensional manifold as $\Lambda$, the operator $\Delta_{\Lambda}$
is equivalent to second derivative in the direction tangent to $\Lambda$, previously
denoted as $d^2_s$.

To establish the preconditioners we follow the framework of operator preconditioning \cite{mardal2011preconditioning}.
That is, having shown stability of \eqref{eq:3D1Dvar} we define the preconditioner as a
Riesz map $\mathcal{B}:X'\rightarrow X$ with respect to inner products/norms in which the problem
was shown to be well-posed. In turn, the framework relates the conditioning of the (preconditioned)
operator $\mathcal{B}\mathcal{A}:X\rightarrow X$ to the stability constants of the Brezzi
theory. More precisely, independence of the constants with respect to a particular problem
parameter translates to robustness of the preconditioner (in the said parameter variations).
Note that in case of \eqref{eq:3D1Dvar}, the critical Brezzi constants are $C_B$ and $\beta$
related to the boundedness and inf-sup condition of the bilinear form $b$, cf. \Cref{th:bnb}.
We remark that to obtain a robust discrete preconditioner, stable discretization is required
in addition to well-posedness of the continuous problem. In particular, the Brezzi constants
of the discretized problem must be independent of the discretization parameter.

Applying operator preconditioning, the analysis of well-posedness in \cite{kuchta2021analysis}
and \Cref{th:bnb} yields to two different preconditioners for \eqref{eq:coupled_3d_1d_operator}.
Here, a key difference is the construction of the multiplier space, as \cite{kuchta2021analysis} consider $H^{-1/2}(\Lambda)$.
However this space does not lead to robust algorithms as we shall demonstrate next. In fact,
intersection spaces of \Cref{th:bnb} will be needed to obtain robustness in material parameters. At the same time,
at the end of this work it will become clear that the radius $\eps$, as a parameter, plays a special and
critical role in the robustness of a preconditioner stemming from \Cref{th:bnb}.


\subsubsection*{\emph{2d-1d} preconditioning example}
\begin{figure}
\begin{center}  
  \includegraphics[height=0.3\textwidth]{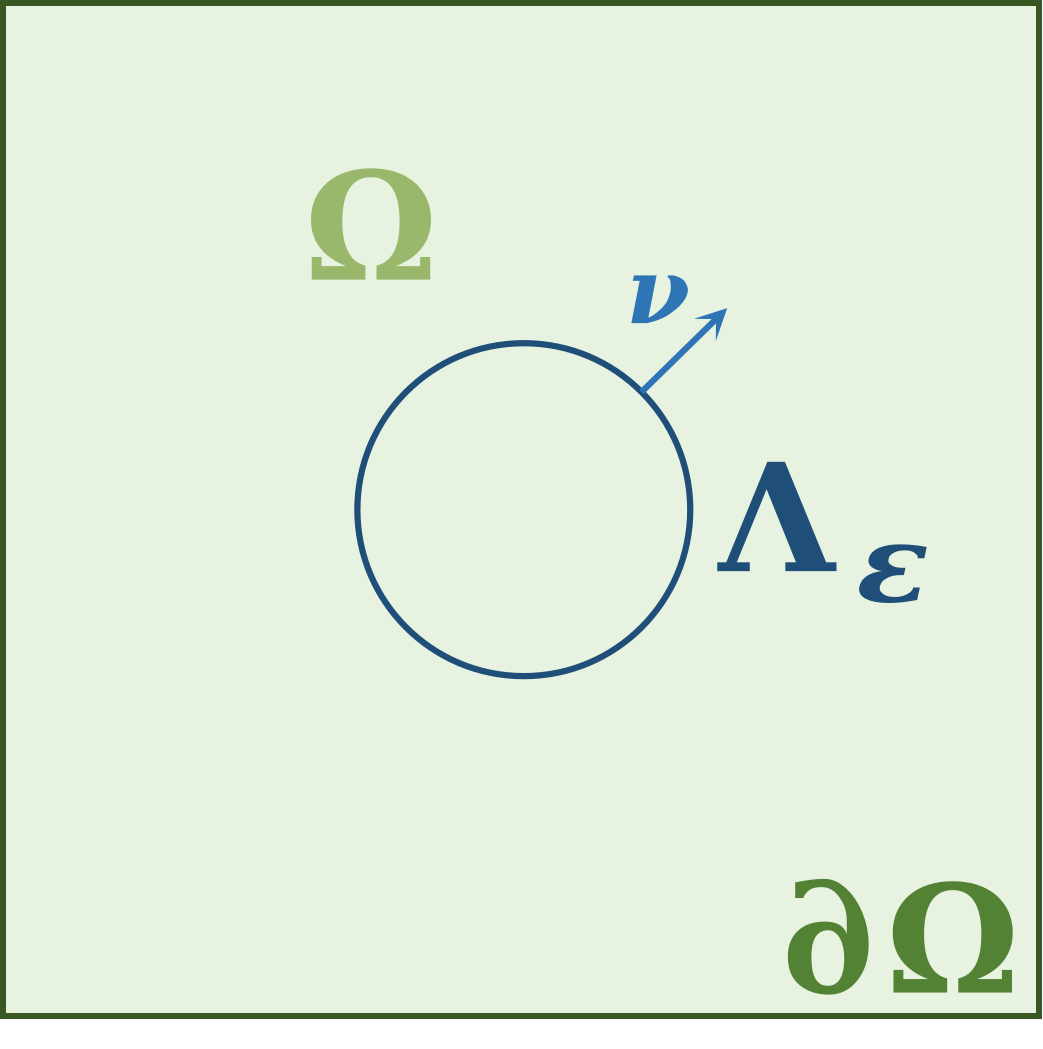}
\hspace{10pt}
  \includegraphics[height=0.3\textwidth]{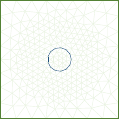}
\vspace{-5pt}
\caption{(Left) Geometrical setting of the \emph{2d-1d} problem \eqref{eq:twoD_oneD_strong}.
  Domain $\Lambda_{\epsilon}$ plays similar role to the centerline in \emph{3d-1d}
  system \eqref{eq:pde}. (Right) In the experiments triangulation of $\Omega$ conforms to $\Lambda$.
}
\label{fig:twoD_oneD_domain}
\end{center}
\end{figure}

Let us illustrate the challenges of developing a robust preconditioner for mixed-dimensional equations
by using a $2d$ analogue of \eqref{eq:coupled_3d_1d_operator} where we let $\Lambda=\Lambda_{\eps}$ be a closed curved
contained in the interior of $\Omega\subset\mathbb{R}^2$, cf. \Cref{fig:twoD_oneD_domain}. 
Here the coupling between the problem unknowns defined on $\Omega$ and $\Lambda$ shall be realized
by the standard $2d$-$1d$ trace operator leading to the coupled problem:
Find $y=[u, u_{\odot}, p_{\odot}]\in Y$ such that
\begin{equation}\label{eq:coupled_2d_1d_operator}
  \mathcal{A}y = f\mbox{ in }Y',\quad
  \mathcal{A}
 = \begin{pmatrix}
  -\Delta & & \tracel'\\
  & -\kappa_{\odot}\Delta_{\Lambda} & -I \\
  \tracel & -I & 0\\
\end{pmatrix}.
\end{equation}
Note that \eqref{eq:coupled_2d_1d_operator} is structurally similar to the
$3d$-$1d$ coupled problem \eqref{eq:coupled_3d_1d_operator}. However, for simplicity,
the model $2d$-$1d$ problem contains only a single parameter, $\kappa_{\odot}>0$, 
whose value can be (arbitrarily) large or small. We also recall that $\Lambda$ is closed and thus, in
order to have an invertible operator on the whole of $H^1(\Lambda)$, we set
$-\Delta_{\Lambda}=-\Delta_{\Lambda}+I$ in this section. We finally remark that in \eqref{eq:coupled_3d_1d_operator}
the Laplacian on $\Lambda$ is invertible due to boundary conditions imposed on $\partial\Omega_{\ominus}$
outside the coupling surface. 

Following \cite{kuchta2016preconditioners} where \eqref{eq:coupled_2d_1d_operator} has been shown
to be well-posed with $Y=H^{1}_0(\Omega)\times H^1(\Lambda)\times H^{-1/2}(\Lambda)$ we shall
construct preconditioners as Riesz maps with respect to two different inner products on $Y$ leading
to operators
\begin{equation}\label{eq:twoD_oneD_preconditioners}
  \mathcal{B}_0 = \begin{pmatrix}
    -\Delta & & \\
    & -\kappa_{\odot}\Delta_{\Lambda} & \\
    & & -\Delta^{-1/2}_{\Lambda}
  \end{pmatrix}^{-1},\quad
  \mathcal{B}_1 = \begin{pmatrix}
    -\Delta & & \\
    & -\kappa_{\odot}\Delta_{\Lambda} & \\
    & & -\Delta^{-1/2}_{\Lambda} - \kappa_{\odot}^{-1}\Delta^{-1}_{\Lambda}
  \end{pmatrix}^{-1}.
  \quad
\end{equation}
Here $\mathcal{B}_0$ stems from the analysis in \cite{kuchta2016preconditioners} and
can be seen to be analogous to the Riesz map preconditioner established in
\cite{kuchta2021analysis} for the coupled $3d$-$1d$ problem
\eqref{eq:coupled_3d_1d_operator}. In particular, the multiplier block of the preconditioner
is a Riesz map of $H^{-1/2}(\Lambda)$ and is independent of the parameter $\kappa_{\odot}$. On the other hand, with $\mathcal{B}_1$ the
multiplier is sought in the intersection space $H^{-1/2}(\Lambda) \cap \kappa_{\odot}^{-1/2}H^{-1}(\Lambda)$, cf.
\Cref{th:bnb}.

We shall compare the two preconditioners in terms of their spectral condition
numbers defined as the ratio between the largest and smallest
in magnitude eigenvalues of the problem: Find $y\in Y$, $\lambda\in\mathbb{R}$
\begin{equation}\label{eq:twoD_oneD_cond}
  \mathcal{A}y
  =
  \lambda
  \mathcal{B}_i^{-1}y\quad\text{ in }Y'.
\end{equation}
To assess robustness of the preconditioners we consider the generalized eigenvalue problem \eqref{eq:twoD_oneD_cond} 
for different values of $10^{-8} \leq \kappa_{\odot} \leq 10^{8}$ and refinements of the domain
$\Omega=(-1, 1)^2$ with $\Gamma_{\eps}=\left\{x\in\Omega, \lvert x\rvert = 0.1 \right\}$.
We then discretize the problem by continuous linear Lagrange elements ($\mathbb{P}_1$) where
the finite element mesh of $\Omega$ always conforms to $\Lambda$, cf. \Cref{fig:twoD_oneD_domain}.

Before addressing the preconditioners let us briefly comment on the approximation property 
of the chosen discretization. To this end we note that \eqref{eq:coupled_2d_1d_operator}
is associated with a system
\begin{equation}\label{eq:twoD_oneD_strong}
\begin{aligned}
  -\Delta u  &= f &\quad\mbox{ in }\Omega,\\
  -\kappa_{\odot}\Delta_{\Lambda} u_{\odot} - \jump{\nabla u}\cdot\nu &= f_{\odot}&\quad\mbox{ on }\Lambda,\\
  u - u_{\odot} &= g &\quad\mbox{ on }\Lambda
\end{aligned}
\end{equation}
which is here supplied with homogeneous Dirichlet boundary conditions for $u$, $u_{\odot}$ on
$\partial\Omega$ and $\partial\Lambda$ respectively. Here $\jump{\cdot}$ is the jump operator on $\Lambda$,
$\jump{v}={v^{+}-v^{-}}$, defined with respect to the normal vector on the curve which points
from the positive to the negative side.

Using \eqref{eq:twoD_oneD_strong}
we measure convergence of the discrete approximations $u_h$, $u_{\odot, h}$, $p_{\odot, h}$
in the norms induced by the (symmetric and positive definite) operator $\mathcal{B}^{-1}_1$. Here, the
linear systems stemming from discretization\footnote{
In all the presented numerical examples FEniCS\cite{LoggMardalEtAl2012}-based module FEniCS\textsubscript{ii}\cite{kuchta2021assembly}
was used to discretize the coupled problems.
} by $\mathbb{P}_1$ elements are  
solved with a preconditioned MinRes solver using relative tolerance
of $10^{-10}$, see also \Cref{fig:twoD_oneD_isect0_isect1}.

\begin{figure}
  \centering
  \includegraphics[height=0.24\textwidth]{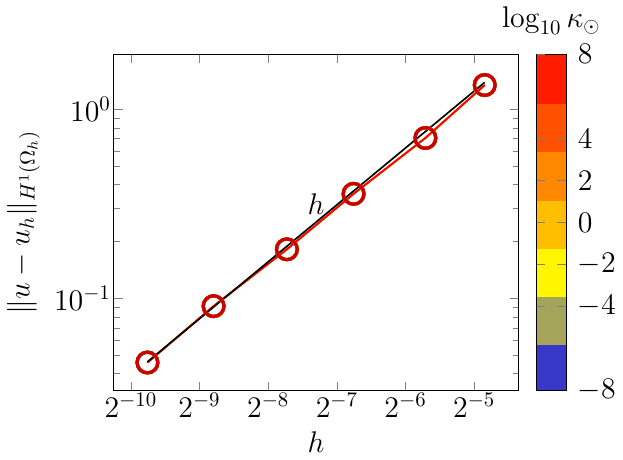}
  \includegraphics[height=0.24\textwidth]{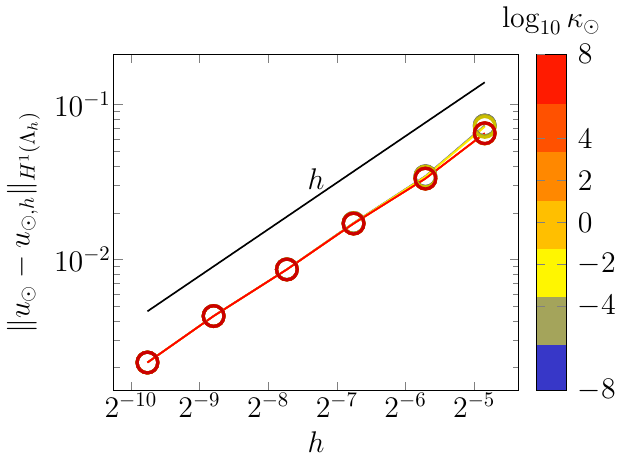}
  \includegraphics[height=0.24\textwidth]{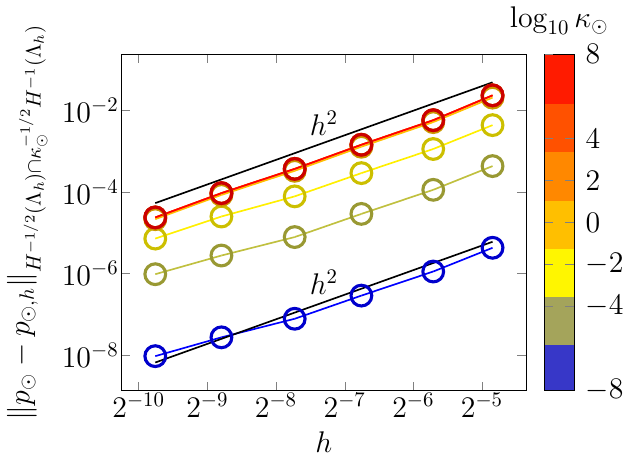}
  \vspace{-20pt}
  \caption{Approximation errors of $\mathbb{P}_1$ discretization of
    \eqref{eq:coupled_2d_1d_operator} for different values of $\kappa_{\odot}$.
  }
  \label{fig:twoD_oneD_cvrg}
\end{figure}

The obtained approximation errors
are plotted in \Cref{fig:twoD_oneD_cvrg}, which, for all the values of $\kappa_{\odot}$, shows
linear convergence in the respective $H^1$-norms for the error in $u$ and $u_{\odot}$. Quadratic convergence
can be seen for the multiplier in the norm of the intersection space
$H^{-1/2}(\Lambda) \cap \kappa_{\odot}^{-1/2}H^{-1}(\Lambda)$. Without including
the results we remark that the intersection norm is essential for obtaining
robust approximation. In particular, we observed that measuring convergence in
only the $H^{-1/2}$-norm yields quadratic convergence for large values of
$\kappa_{\odot}$ while for small values the rate drops to linear.

\begin{figure}
  \centering
  \includegraphics[height=0.35\textwidth]{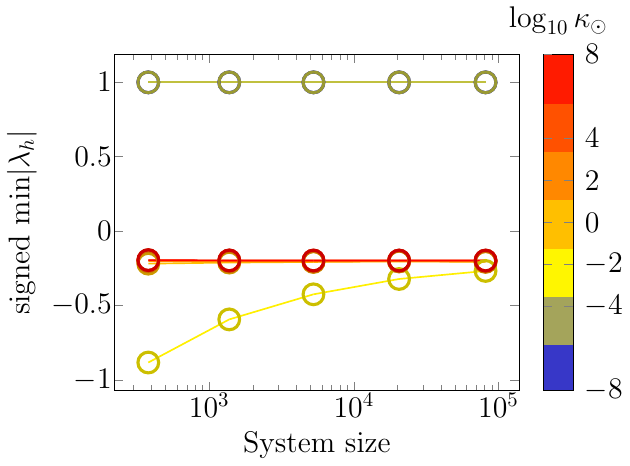}
  \hspace{15pt}
  \includegraphics[height=0.35\textwidth]{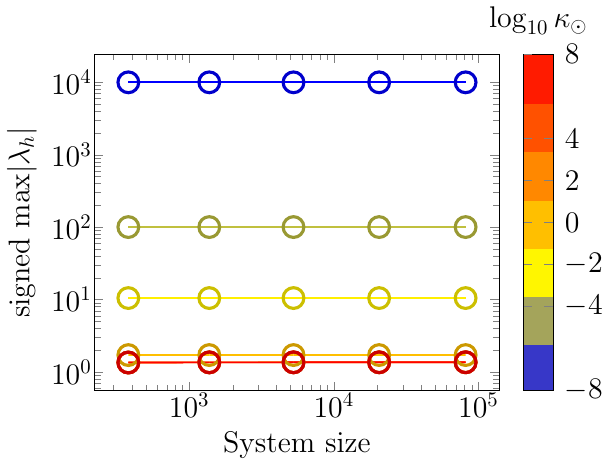}
  \vspace{-5pt}
  \caption{
    Performance of preconditioner $\mathcal{B}_0$ from \eqref{eq:twoD_oneD_preconditioners}
    for problem \eqref{eq:coupled_2d_1d_operator} with varying $\kappa_{\odot}$.
    Here $\Omega=(-1, 1)^2$, $\Gamma_{\eps}$ is a circle of radius $\eps=0.1$.
    All the finite element spaces are discretized by $\mathbb{P}_1$ elements. The preconditioner
    does not yield $\kappa_{\odot}$-bounded condition numbers.
  }
  \label{fig:twoD_oneD_isect0}
\end{figure}

\begin{figure}
  \centering
  \includegraphics[height=0.35\textwidth]{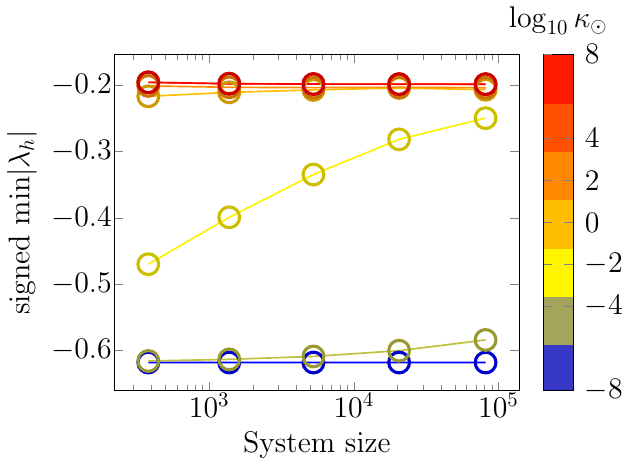}
  \hspace{15pt}
  \includegraphics[height=0.35\textwidth]{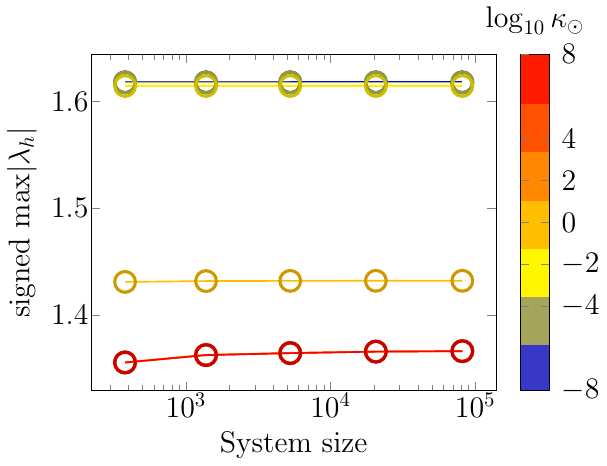}
  \vspace{-5pt}  
  \caption{
    Performance of preconditioner $\mathcal{B}_1$ from \eqref{eq:twoD_oneD_preconditioners}
    for $2d$-$1d$ coupled problem \eqref{eq:coupled_2d_1d_operator} with varying $\kappa_{\odot}$. Setup as
    in \Cref{fig:twoD_oneD_isect0}. Preconditioner is robust in $\kappa_{\odot}$ and mesh size.
  }
  \label{fig:twoD_oneD_isect1}
\end{figure}

\begin{figure}
  \centering
  \includegraphics[height=0.35\textwidth]{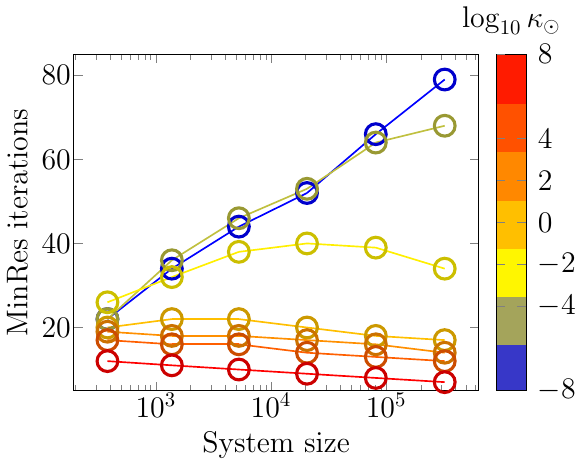}
  \hspace{15pt}
  \includegraphics[height=0.35\textwidth]{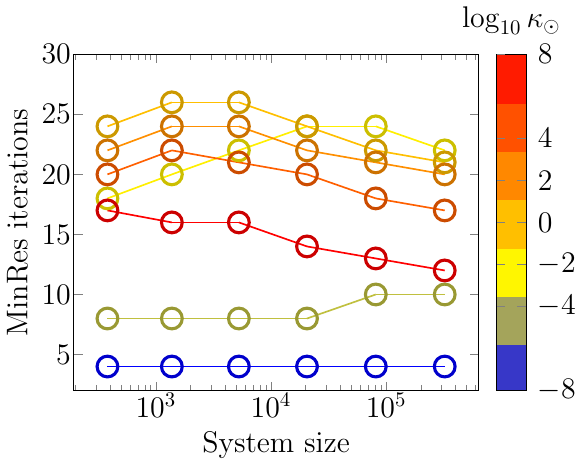}
  \vspace{-5pt}  
  \caption{
    Convergence of \eqref{eq:twoD_oneD_preconditioners}-preconditioned MinRes solver
    for problem \eqref{eq:coupled_2d_1d_operator} with varying $\kappa_{\odot}$.
    Only the preconditioner $\mathcal{B}_1$ (right) reflecting $p_{\odot}\in H^{-1/2}\cap\kappa_{\odot}^{-1/2}H^{-1}$
    is parameter robust. Setup as in \Cref{fig:twoD_oneD_isect0}.
  }
  \label{fig:twoD_oneD_isect0_isect1}
\end{figure}

Having verified the discretization scheme (and our implementation) we return to
preconditioning and stability of the eigenvalue problem \eqref{eq:twoD_oneD_cond}.
We summarize the results in \Cref{fig:twoD_oneD_isect0} and \Cref{fig:twoD_oneD_isect1}
which show the extrema, i.e. minimum and maximum absolute values, of the eigenvalues $\lambda_h$
of the discretized problem \eqref{eq:twoD_oneD_cond}. We note that the bounds are
plotted together with their sign which carries relevant information related e.g.
to the discrete inf-sup or coercivity conditions.

Using $\mathcal{B}_0$, the extremal eigenvalues of \eqref{eq:twoD_oneD_cond} are
displayed in \Cref{fig:twoD_oneD_isect0}. We observe that for each $\kappa_{\odot}$ the quantities
are bounded in $h$ verifying stability of the discrete problem with $\mathbb{P}_1$ elements
and the norm induced on $Y$ by $\mathcal{B}^{-1}_0$, see \cite{kuchta2016preconditioners}.
However, there is an apparent growth of the largest in magnitude eigenvalue as
$\kappa_{\odot}$ becomes small. As such $\mathcal{B}_0$ does not yield
parameter robust solver.
On the other hand, the preconditioner $\mathcal{B}_1$ yields spectral bounds that
are stable in $\kappa_{\odot}$ and mesh refinement. We note that for small $\kappa_{\odot}$ the
observed upper and lower bounds are close to their theoretical values \cite{murphy2000note, rusten1992preconditioned}
of $(1+\sqrt 5)/2$ and $(1-\sqrt 5)/2$ respectively. This is due to the multiplier norm being
then dominated by the matrix of the $H^{-1}$-inner product such that the discrete
preconditioner (which computes $-\Delta_{\Lambda}^{-1}$ by LU decomposition) is close to being the
exact Schur complement preconditioner.

To illustrate how the conditioning translates to performance of iterative
solvers \Cref{fig:twoD_oneD_isect0_isect1} reports the iteration counts
of the MinRes solver using the two preconditioners \eqref{eq:twoD_oneD_preconditioners}.
Here we reuse the setup of the previous eigenvalue experiments while the right
hand side $f$ in \eqref{eq:coupled_2d_1d_operator} is based on the manufactured
solution setup described above. For each value of $\kappa_{\odot}$ and the mesh size the solver
is started from 0 initial guess and terminates once the preconditioned residual norm
is reduced by factor $10^{10}$. Both preconditioners are computed exactly using
LU for the two leading blocks while the multiplier block, in particular the fractional
term $-\Delta_{\Lambda}^{-1/2}$ is realized via spectral decomposition\footnote{
Spectral decomposition is not suitable for practical applications because of its cubic
scaling. However, Riesz maps in intersections of fractional order Sobolev spaces
can be approximated with optimal complexity by rational approximations \cite{budivsa2022rational}.
}, see \cite{kuchta2016preconditioners}
for precise definition.

In \Cref{fig:twoD_oneD_isect0_isect1} it can be seen
that for large values of $\kappa_{\odot}$ both preconditioners yield iterations which are stable in mesh refinement.
However, for small values, i.e. when the term $\kappa_{\odot}^{-1}\Delta^{-1}_{\Lambda}$ becomes large,
$\mathcal{B}_0$ shows dependence on the parameter. We conclude that the intersection space exploited in the definition of $\mathcal{B}_1$ is crucial for parameter robustness.


\section{Definition of a preconditioner for the \emph{3d-1d} problem: performances and drawbacks}\label{sec:3d1d_precond}

At this point let us return to the original coupled $3d$-$1d$ problem \eqref{eq:coupled_3d_1d_operator}
that we shall now consider with a preconditioner
\begin{equation}\label{eq:coupled_3d_1d_preconditioner}
  \mathcal{B} = \begin{pmatrix}
    -\kappa\Delta & & \\
    & -\kappa_{\odot}\eps^2\Delta_{\Lambda} & \\
    & & -\frac{\eps^2}{\kappa}\Delta^{-1/2}_{\Lambda} - \kappa_{\odot}^{-1}\Delta^{-1}_{\Lambda}
  \end{pmatrix}^{-1},
\end{equation}
i.e. the Riesz operator associated to the inner product of the space
in which well-posedness of \eqref{eq:coupled_3d_1d_operator} was shown
in \Cref{th:bnb}.

\begin{figure}
  \centering
  \includegraphics[height=0.42\textwidth]{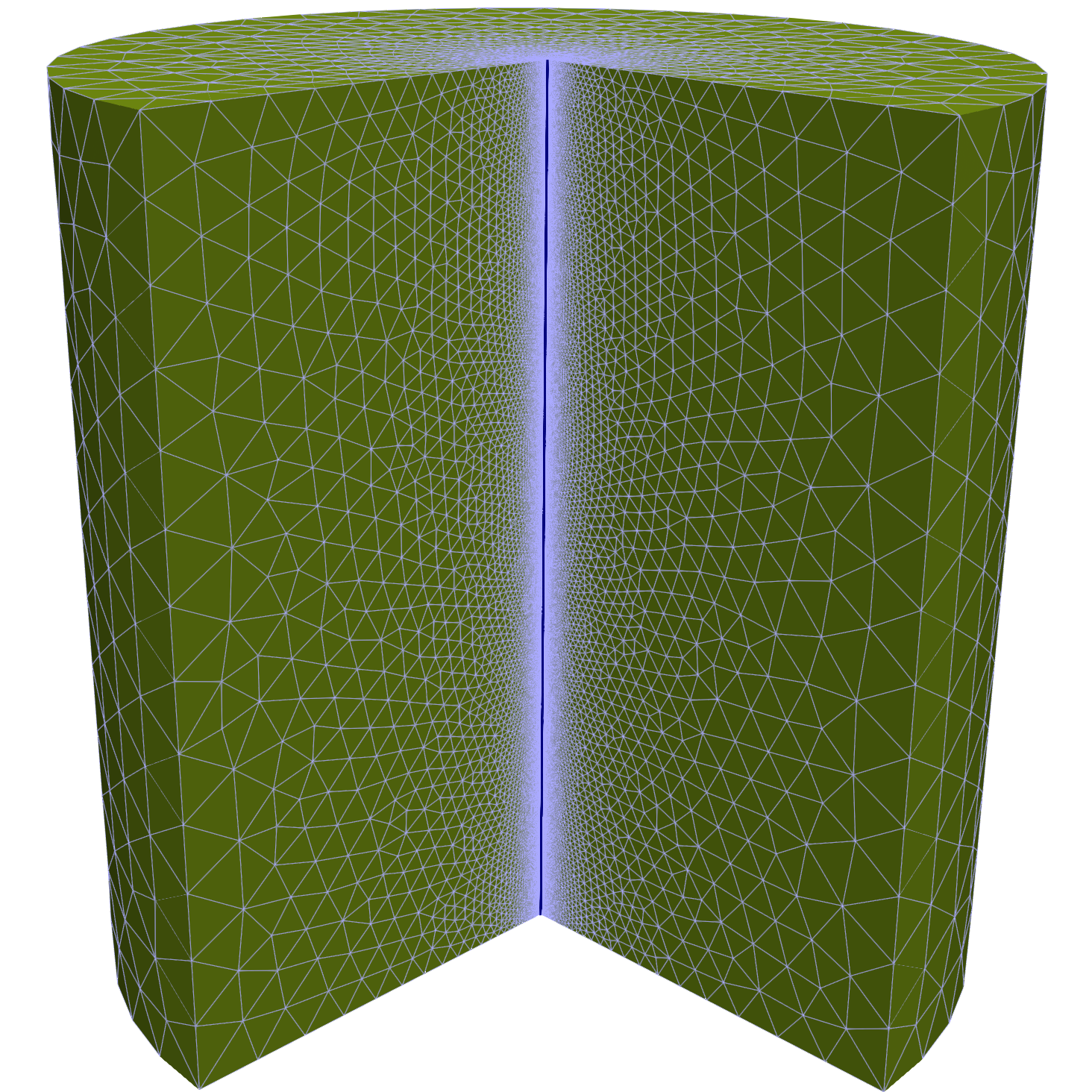}
  \hspace{15pt}
  \includegraphics[height=0.42\textwidth]{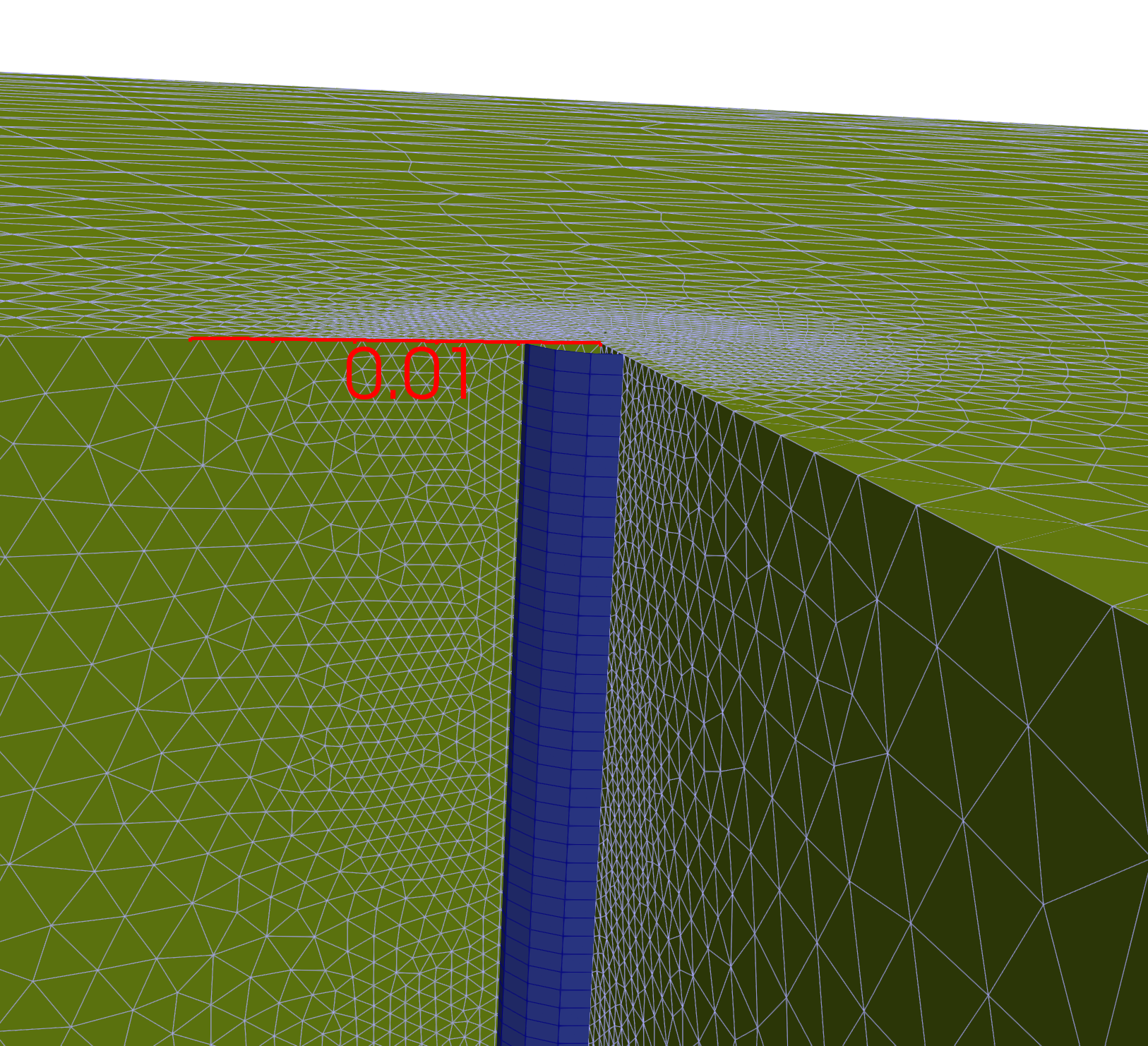}
  \vspace{-5pt}  
  \caption{
    Computational mesh considered in experiments for the coupled
    $3d$-$1d$ problem \eqref{eq:coupled_3d_1d_preconditioner}. Here $\eps=2\cdot10^{-3}$.
    Mesh of $\Omega$ conforms to the centerline $\Lambda$ as well as to the coupling surface
    $\Gamma_{\eps}$ leading to high refinement near $\Lambda$ (edges of the
    mesh elements are pictured in light blue color). (Right) Virtual surface used for computing $\mtrace$ is rendered
    in dark blue color.    
  }
  \label{fig:threeD_oneD_mesh}
\end{figure}

In order to test \eqref{eq:coupled_3d_1d_preconditioner} we consider cylindrical
domains $\Omega_{\oplus}$, $\Omega_{\ominus}$ each with height 1 and radii of 0.5 and $\eps$
respectively. Following \cite{kuchta2021analysis} we discretize $\Omega=\Omega_{\oplus}\cup\Omega_{\ominus}$
such that the mesh conforms both to the interface $\Gamma_{\eps}$ and the centerline $\Lambda$,
see \Cref{fig:threeD_oneD_mesh}. We remark that the conformity assumption was used in \cite{kuchta2021analysis}
to show stability of the discrete problem with $\mathbb{P}_1$ elements (used below). At the same time,
the assumption leads to greatly refined meshes in the vicinity of $\Gamma_\eps$. It also
increases the cost of mesh generation and, due to
the size of the resulting system, restricts\footnote{
As an example, for radius $\eps=1\cdot10^{-2}$ the finest mesh considered
contained roughly 11 million tetrahedra. Using $\mathbb{P}_1$ elements the
number of $3d$ unknowns is then $\sim\!\!2$ million while the multiplier space has
$\sim\!\!2$ thousand degrees of freedom.
} the type of experiments we can perform (on our
serial computational setup). In the following we shall thus limit
the computational study only to robustness of iterative methods. We remark that
stable discretization of \eqref{eq:coupled_3d_1d_preconditioner} is possible
also if the mesh of $\Omega$ is independent of $\Lambda$ and $\Gamma$, cf. \cite{kuchta2021analysis}.
However, we argue that the conforming setup is simpler and more transparent.

\begin{figure}
  \centering
  \includegraphics[height=0.35\textwidth]{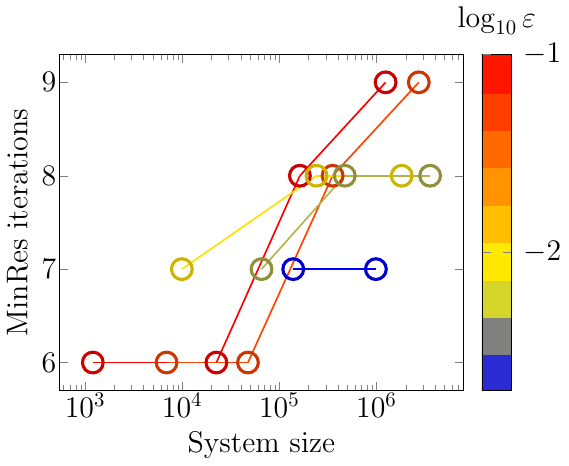}
  \hspace{15pt}
  \includegraphics[height=0.35\textwidth]{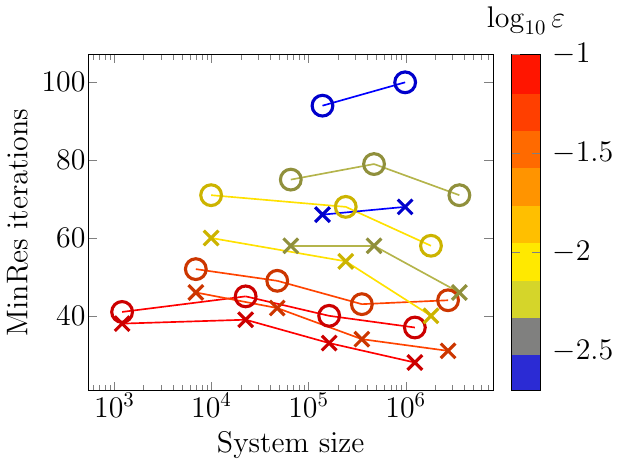}
  \vspace{-5pt}  
  \caption{
    Performance of preconditioner \eqref{eq:coupled_3d_1d_preconditioner} for the
    coupled $3d$-$1d$ problem \eqref{eq:coupled_3d_1d_operator} and varying radius.
    (Left) We set $\kappa=10^{8}$, $\kappa_{\odot}=10^{-8}$. (Right) We set
    $\kappa=10^{-8}$, $\kappa_{\odot}=10^{8}$. The leading block of the preconditioner
    uses single AMG V-cycle ($\circ$ markers) or 10 PCG iterations with AMG preconditioner ($\times$ markers).
    The remaining blocks are realized by LU.
  }
  \label{fig:threeD_oneD_precond}
\end{figure}

In \Cref{fig:threeD_oneD_precond} we report on the convergence of the MinRes solver
using preconditioner \eqref{eq:coupled_3d_1d_preconditioner} for $2\cdot10^{-3}\leq \eps \leq 10^{-1}$
and two parameter regimes. Here the action of the leading block $-\kappa\Delta^{-1}$
of \eqref{eq:coupled_3d_1d_preconditioner} is approximated 
in terms of (i) a single V-cycle of algebraic multigrid (AMG)
and (ii) 10 steps of preconditioned conjugate gradient (PCG) method using AMG as preconditioner.
The remaining two blocks of \eqref{eq:coupled_3d_1d_preconditioner} are
computed by LU factorization. We remark that the choice (i) is more
practical while with (ii) the preconditioner is almost exact as the absolute
residual norm after the 10 PCG steps is typically $<10^{-15}$ in our case. With
(ii) we thus aim to ensure that the effects of parameter variations on MinRes
convergence are (mostly) due to the construction of the multiplier preconditioner.
In both cases the convergence criterion for the MinRes solver requires reducing
the preconditioned residual norm by a factor $10^{10}$. Finally, the coupling operator
$\mtrace$ is approximated using a Legendre quadrature of degree 20.

Considering the results in \Cref{fig:threeD_oneD_precond} we observe that
performance of the preconditioner differs dramatically between the two regimes. 
When $\kappa=10^{8}$, $\kappa_{\odot}=10^{-8}$, such that
the $H^{-1}$-term can dominate the multiplier block in $\mathcal{B}$, the iteration
counts are practically independent of $\eps$. Note that here only the construction
(i) for the leading block is considered as it already yields low enough iterations. On the other hand, with
$\kappa=10^{-8}$, $\kappa_{\odot}=10^{8}$ the solver performance
deteriorates for small radii. This is true for the construction (i) which,
for the different radii $\eps\in\left\{1\cdot10^{-1}, 5\cdot10^{-2}, 1\cdot10^{-2}, 5\cdot10^{-3}, 2\cdot10^{-3}\right\}$
and the finest refinement levels yields the iterations counts
of 37, 44, 58, 71, 100 as well as for (ii) where convergence is reached respectively
after 28, 31, 40, 46, 68 iterations.

The unbounded iterations, in particular with
preconditioner (ii), bring into question the stability of the coupled $3d$-$1d$ problem
\eqref{eq:coupled_3d_1d_operator} with preconditioner \eqref{eq:coupled_3d_1d_preconditioner} and
in particular the intersection space $Q_{\odot}$ for the multiplier.
The lack of robustness is surprising as it suggests that
radius $\eps$ in \eqref{eq:coupled_3d_1d_operator} does not behave as a standard
material parameter in the sense that corresponding weighting in the (appropriate) intersection
space does not yield robustness with respect to its variations.
However, this result may be qualitatively justified on the basis of the stability analysis,
where we noticed the influence of the constant $C_{IT}$ in the inf-sup condition. We will
show in the next section that $C_{IT}$ may depend on $\eps$, which may explain the
findings in \Cref{fig:threeD_oneD_precond}.

\subsubsection*{Varying radius in \emph{2d-1d} preconditioning example}

In order to gain more insight into the observation that in the $3d$-$1d$ setting of \eqref{eq:coupled_3d_1d_operator}
the multiplier preconditioner reflecting the intersection space $\kappa^{-1/2}\eps H^{-1/2}\cap \kappa^{-1/2}_{\odot}H^{-1}$
did not lead to robustness, let us return to the $2d$-$1d$ operator \eqref{eq:coupled_2d_1d_operator}.
Next, we shall investigate the properties of the preconditioner $\mathcal{B}_1$
when $\eps$ varies. We recall that $\mathcal{B}_1$ was found to be robust with respect
to the material parameters.

\begin{figure}
  \centering
  \includegraphics[height=0.245\textwidth]{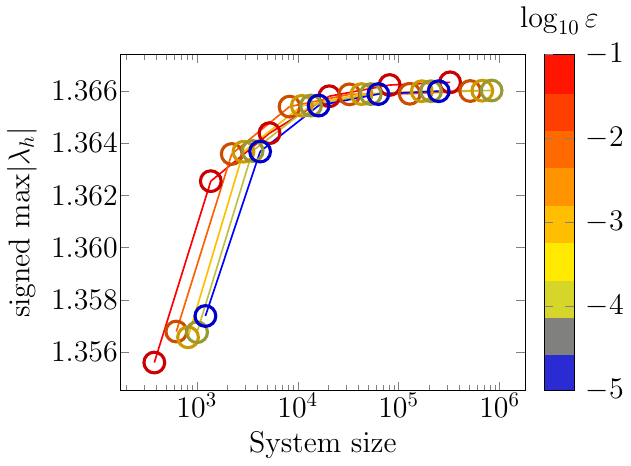}
  \includegraphics[height=0.245\textwidth]{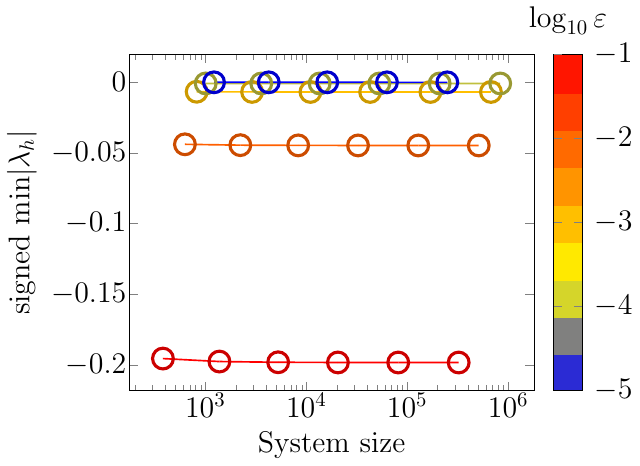}
  \includegraphics[height=0.245\textwidth]{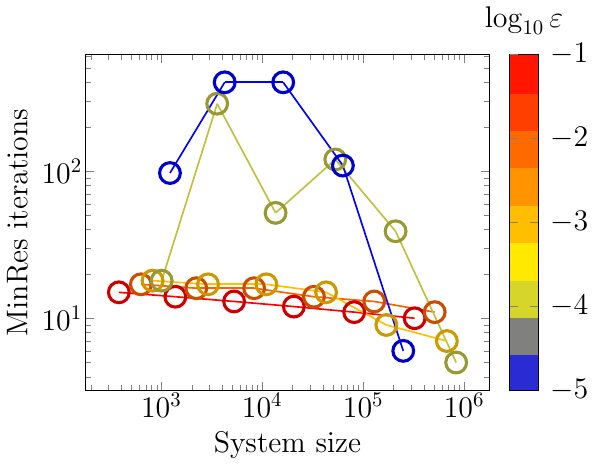}
  \vspace{-20pt}
  \caption{
  Performance of preconditioner $\mathcal{B}_1$ in \eqref{eq:twoD_oneD_preconditioners}
  for problem \eqref{eq:coupled_2d_1d_operator} with $\kappa_{\odot}=10^{10}$,
  $\Omega=(-1, 1)^2$, $\Gamma_{\eps}$ is a circle of radius $10^{-5} \leq \eps \leq 10^{-1}$.
  All the finite element spaces are discretized by $\mathbb{P}_1$ elements.
  (Left, center) The conditioning 
  deteriorates with $\eps$ as the largest in magnitude eigenvalue \eqref{eq:twoD_oneD_cond} is bounded in the
  parameter while the smallest in magnitude eigenvalue decreases with $\eps$.
  (Right) For small radii the iterations are not stable (maximum number of allowed iterations is
  set to 400). 
  }
  \label{fig:twoD_oneD_isect1_radius}
\end{figure}

Using the previous experimental setup, \Cref{fig:twoD_oneD_isect1_radius}
shows the performance of the preconditioner when $\kappa_{\odot}=10^{10}$ and $\eps\leq 10^{-1}$.
Here, by the choice of a large $\kappa_{\odot}$ we wish to put emphasis on the fractional
part of the multiplier preconditioner, which, analogously to \Cref{th:bnb}, brings
the inverse trace constant into the Brezzi estimates.
In \Cref{fig:twoD_oneD_isect1_radius} we observe that while the largest eigenvalues
are bounded in $\eps$, the smallest in magnitude eigenvalue approaches 0 as
$\Lambda_{\eps}$ shrinks. We note that for all $\eps$ the two eigenvalue bounds 
are stable in mesh size $h$.
In the figure we further report convergence history of the MinRes solver (run with the
same settings as in the previous experiments).
We observe that for small radii the iteration counts grow rapidly with initial mesh
refinement before decreasing back on finest meshes, cf. bounded iterations
with mesh refinement in \Cref{fig:twoD_oneD_isect0_isect1}. However, the
limit $\eps\rightarrow 0$ does not seem to affect the iterative solver as clearly
as the blow up of the condition number, cf. $\text{min}\lvert \lambda_h \rvert$ in \Cref{fig:twoD_oneD_isect1_radius}.
We attribute this behaviour to the specific choice of the right-hand side and 0 initial
guess in our numerical setup. Nevertheless, the sensitivity to the radius furthers our claim
of the special role of $\eps$ in understanding preconditioner robustness.




\section{The role of the inner radius on mixed-dimensional problems}\label{sec:3d1d_issue}

As discussed in the preceding sections, it has been observed that the robustness of the preconditioner \eqref{eq:coupled_3d_1d_preconditioner} diminishes as the inner radius of the inclusion approaches zero. This behavior is not limited to the $3d$-$1d$ preconditioner formulation but is also evident in the $2d$-$1d$ examples. Consequently, it suggests that the detrimental effect of $\eps$ is not a result of the dimensional reduction technique itself but rather inherent in the mathematical structure of the problem.

Specifically, the dependence on $\eps$ is manifested in the constant of the inf-sup stability property. The analysis in Theorem \ref{th:bnb} reveals that the boundedness constant $C_B$ of the bilinear form $b$ remains independent of the inner diameter. However, the inf-sup constant $\beta$ is not guaranteed to be independent of $\eps$ as the presence of the inverse trace constant $C_{IT}$ introduces questions regarding its independence from the inner radius. The inverse trace theorem, a well-established result in functional analysis, provides insights into this matter (see, for example, \cite{steinbach2007numerical}, \cite{necas2011direct}, \cite{lions2012non}). However, deriving an accurate bound for $C_{IT}$ is a challenging task, particularly when considering a parameterized domain.

In this section, our objective is to investigate the relationship between the inf-sup constant $\beta$ and the parameter $\eps$, which is closely associated with the inverse trace constant $C_{IT}$. To accomplish this, we analyze a simplified yet significant scenario that allows us to elucidate the connection between the inf-sup constant and the inner radius.

\subsection{The $2d$-$1d$ formulation for the perforated domain problem}
  \begin{figure}[b]
  \centering
   \hspace{15pt}
  \includegraphics[width=0.8\textwidth]{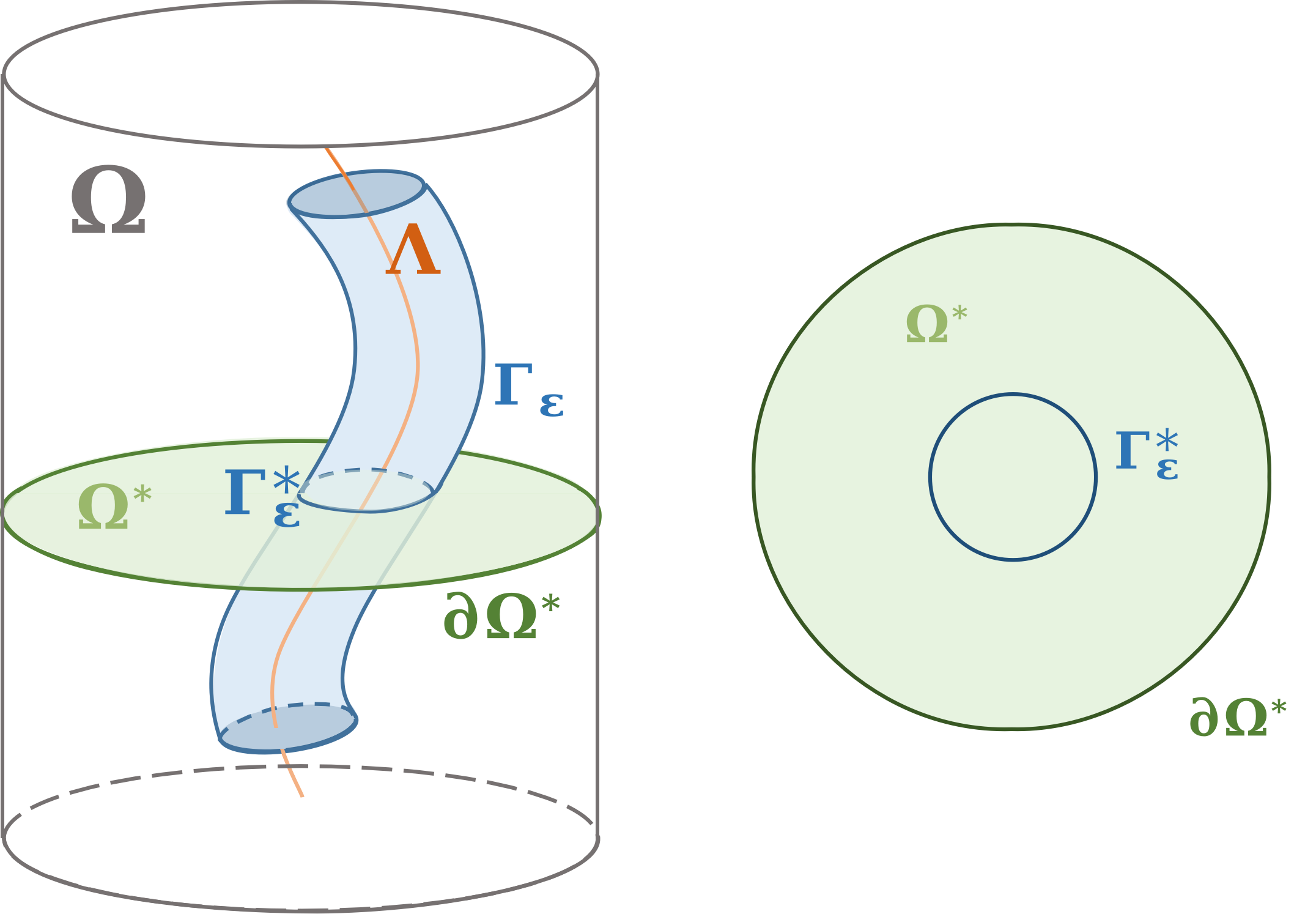}
\hspace{50pt}
\vspace{-5pt}
  \caption{
    (Left) Pictorial representation of a vessel with radius $\eps$, centerline $\Lambda$ and boundary
    $\Gamma_\eps$ immersed in a three-dimensional domain $\Omega$. By cutting a slice of $\Omega$, we
    obtain the two-dimensional domain $\Omega^*$ with boundary $\partial \Omega^*$ and circular one-dimensional
    inclusion $\Gamma^{*}_\eps$.   (Right) Pictorial representation of the slice-domain $\Omega^*$.
  }
  \label{fig:cont_u}
\end{figure}

  Let us consider a generalized cylindrical vessel immersed in a three-dimensional domain $\Omega$. The
  vessel surface is denoted by $\Gamma_{\eps}$, where $\eps$ indicates the cylinder radius.
  We are interested in studying how the mathematical structure of the problem behaves when $\eps \rightarrow 0$.
  A straightforward approach can be to select a slice of $\Omega$, here denoted by the superscript "$*$", 
  such that we obtain a domain $\Omega^*$, where the curve $\Gamma_{\eps}^{*}$ is the restriction of $\Gamma_{\eps}$ on the slice, see \Cref{fig:cont_u}. Clearly, $\text{diam}\Gamma_\eps^*$ depends on $\eps$. In turn, studying the effects of shrinking
  $\text{diam}\Gamma_\eps^*$ in $\Omega^*$ (a $2d$-$1d$ system) is representative of the effects of the diminishing vessel/inner radius. The underlying assumption is the complete decoupling of the radius influence
  from the axial direction, which seems reasonable when considering cylindrical setting.
  Indeed, inspecting the expression of the Laplacian in cylindrical coordinates
\[
    \Delta  = \frac{1}{\rho}\frac{\partial}{\partial \rho}\left(\rho \frac{\partial}{\partial \rho} \right) + \frac{1}{\rho^2}\frac{\partial^2}{\partial \phi^2} + \frac{\partial^2}{\partial z^2}
\]
corroborates the assumption: $\rho$ (and in turn $\eps$) does not appear in the axial derivative part. From now, until the end of this section, the superscript "$*$" will be omitted, so that $\Omega$ and $\partial \Omega$ will refer
respectively to the $2d$ sliced domain and its outer boundary. Moreover, $\Gamma_{\eps}$ will
denote the one-dimensional closed curve embedded in $\Omega^*$.  

Let us consider the following $\eps$-dependent Poisson problem defined on $\Omega$: 

\begin{equation}
\begin{aligned}
-\Delta u&=0 & \text { on } \Omega, \\
\left.u\right|_{\Gamma_{\eps}}&=g & \text { on } \Gamma_{\eps}, \\
u&=0 & \text { on } \partial \Omega, \\
\end{aligned}
\label{eq:cont_u}
\end{equation}
where $\Gamma_{\eps}\cap\partial\Omega=\emptyset$ and $u$ is a two-dimensional scalar field defined on $\Omega$.
Now consider a  mixed weak formulation of \eqref{eq:cont_u}, where the 
boundary condition on $\Gamma_{\eps}$ is enforced weakly by a Lagrange multiplier
in order for the the mathematical structure of \eqref{eq:cont_u} to mirror the $3d$-$1d$ problem.
The variational problem reads: 
\begin{equation}
\begin{aligned}
\left(\nabla u, \nabla v\right)_{\Omega}+ (p_\odot, \mathcal{T}v )_{\Gamma_{\eps}} & = 0 & \quad \forall u, v \in H_0^{1}(\Omega) \\
( \mathcal{T}u, q_\odot)_{\Gamma_{\eps}} & = ( g, q_\odot)_{\Gamma_{\eps}} & \quad \forall p_\odot, q_\odot \in H^{-1 / 2}(\Gamma_{\eps})
\end{aligned}.
\label{eq:weak_P}
\end{equation}
Note that the solution operator $\mathcal{A}$ of \eqref{eq:weak_P} has the following block structure
\[
\mathcal{A} = \begin{pmatrix}
  A & B'\\
  B & 0
\end{pmatrix},\quad
\langle Au, v \rangle_{\Omega} = \int_{\Omega}\nabla u\cdot \nabla v,\quad
\langle Bu, q_\odot \rangle_{\Gamma_{\eps}} = \int_{\Gamma_{\eps}} p_\odot v.
\]
%

The question we address next is that of well-posedness of the variational formulation \eqref{eq:weak_P} when $\eps \rightarrow 0$.
In the framework of saddle-point problems, the boundedness of the bilinear
forms $a, b$ in case of \eqref{eq:weak_P} can be easily established with the respective constants
equal to 1. Regarding the existence of the Brezzi inf-sup constant $\beta$,
a straightforward proof will make use of the following theorem \cite{steinbach2007numerical}:
\begin{theorem}\textbf{(Inverse Trace Theorem)} 
Given the necessary smoothness assumption for $\Omega$, the trace operator $\mathcal{T}: H^1(\Omega) \rightarrow H^{1 / 2}(\Gamma_{\eps})$ has a continuous right inverse operator
$$
\mathcal{E}: H^{1 / 2}(\Gamma_{\eps}) \rightarrow H^1(\Omega)
$$
satisfying $\mathcal{T} \mathcal{E} w_\odot=w_\odot$ for all $w_\odot \in H^{1 / 2}(\Gamma_{\eps})$ as well as
$$
\|\mathcal{E} w_\odot\|_{H^1(\Omega)} \leq C_{I T}\|w_\odot\|_{H^{1 / 2}(\Gamma_{\eps})} \quad \forall w_\odot \in H^{1 / 2}(\Gamma_{\eps}).
$$
\end{theorem}
Then, by taking $v=\mathcal{E}v_{q_{\odot}}$, where, by $v_{q_{\odot}}$, we intend an element of $H^{1/2}(\Gamma_{\eps})$ such that the following Riesz mapping properties hold

\begin{equation*}
    \langle v_{q_{\odot}}, w_\odot\rangle_{H^{1/2}(\Gamma_{\eps})} = (q_\odot, w_\odot)_{\Gamma_{\eps}}\text{ and }
    \| v_{q_{\odot}} \|_{H^{1/2}(\Gamma_{\eps})} =\| q_\odot \|_{H^{-1/2}(\Gamma_{\eps})},
\end{equation*}
we obtain that 
\begin{equation*}
\begin{aligned}
\sup _{v \in H^1_0(\Omega), v \neq 0} \frac{( \mathcal{T}v, q_\odot )_{\Gamma_{\eps}}}{\|v\|_{H^1_0(\Omega)}} &  \geq \frac{( v_{q_{\odot}}, q_\odot )_{\Gamma_{\eps}}}{\|\mathcal{E}v_{q_{\odot}}\|_{H^1_0(\Omega)}}
 = \frac{\langle v_{q_{\odot}}, v_{q_{\odot}} \rangle_{H^{1/2}(\Gamma_{\eps})}}{\|\mathcal{E}v_{q_{\odot}}\|_{H^1_0(\Omega)}} \geq \frac{1}{C_{IT}}\| v_{q_{\odot}}\|_{H^{1/2}(\Gamma_{\eps})}.
\end{aligned}
\end{equation*}
What is apparent form the above calculation, is, again, the structural relationship between the existence of the inf-sup constant $\beta$, and the inverse trace inequality constant $C_{IT}$ i.e $\beta=1/C_{IT}$. This suggests (and it is what we want to uncover) a tight relationship between $\beta$ and the inner radius by means of the \textit{trace inequality constant} $C_T$.
Indeed from  Lemma 2.2 in \cite{KVWZ}, the $\eps$ influence is made clear:

\begin{lemma}[Lemma 2.2 in \cite{KVWZ}]
 Let $B_\eps \subset \mathbb{R}^2$ be a circle with a sufficiently small radius $\eps$ and $v \in H^1_0(\Omega)$. Then we have:
\begin{equation}\label{eq:trace_constant}
    \left\| \TT v \right\|_{L^2(\partial B_\eps)} \leq C_T\left\|v\right\|_{H^1_0(\Omega)} = C \sqrt{\eps|\log \eps|} \left\|v\right\|_{H^1_0(\Omega)} 
\end{equation}
\textit{with $C$ positive and independent of $\eps$.}
\end{lemma}

By considering that the extension operator, $\mathcal{E}$, is the right inverse of the trace, $\TT$, the following claim is proposed: \textit{the vessel radius $\eps$ directly affects the inf-sup constant $\beta$ through the trace inequality constant $C_T$}. This claim sounds reasonable if one considers that $\beta$ directly characterizes the lower bound of the bilinear form $b$,
which depends strongly on the trace operator; moreover, it is corroborated by the fact that $\beta= 1/C_{IT}$ so that, assuming a relation between the trace inequality and the inverse trace one, we get that $\eps$ directly affects $\beta$.
That said, in order to prove the claim, it would be required to track in detail the
exact expression for $C_{IT}$, obtaining, thus, a direct link between $\beta$
and $\eps$. Establishing such an analytical expression for $C_{IT}$ can be a
very intricate task, especially when the domain size needs to be considered as the
parameter. Nonetheless, we can rely on numerical analysis and employ
the computability of the constants in a discretized setting, as follows.\\ \\
\textbf{\textit{Brezzi inf-sup constant evaluation.}} Let $V=H^1_0(\Omega)$, $Q=H^{-1/2}(\Gamma_{\eps})$
and $\lVert \cdot \rVert_V$ be the $H^1_0$-norm (induced by the operator $A$) while $\|\cdot\|_Q$ shall be the
$H^{-1/2}_{\Gamma}$-norm induced by the fractional operator $N:=(-\Delta+I)^{-1/2}$ (where both the terms are understood to be defined on $\Gamma_{\eps}$). Consider now a family  $(V_h,\,Q_h)$ of discretization of $(V,\,Q)$ characterized by the following discrete Brezzi inf-sup constant 
\begin{equation}
\sup _{v \in V_h, v \neq 0} \frac{b(v, q_\odot)}{\|v\|_{V_h}} \geq \beta_h\|q_\odot\|_{Q_h} \quad \forall q_\odot \in Q_h.
\label{eq:discrete-inf-sup}
\end{equation}
If $(V_h,\,Q_h)$ is a \textit{stable} discretization i.e. both the discrete Brezzi inf-sup and coercivity constants are bounded from below by a constant which is independent from $h$ \cite{arnold2009stability}, \cite{rognes2012automated}:
\begin{equation}
    \{\beta_h\}_{h \rightarrow 0} \geq \beta > 0,
\end{equation} 
then an approximation of $\beta$ can be obtained by considering the truncated limit $\{\beta_h\}_{h\rightarrow\delta}$, for $\delta$ sufficiently small. 
With the purpose of computing $\beta_h$, we recall the following lemma by Qin \cite{qin1994convergence}:

\begin{lemma}[Qin \cite{qin1994convergence}]
Given a stable discretization $(V_h,\,Q_h)$ for the saddle point problem \eqref{eq:weak_P}, consider the following generalized eigenproblem: Find  $\lambda \in \R$, $0 \neq (u,\,p_\odot) \in V_h \times Q_h$ such that
\begin{equation}
\left\langle u_h, v\right\rangle_V+b\left(v, p_{\odot}\right)+b\left(u_h, q_\odot\right)=-\lambda\left\langle p_{\odot}, q_\odot \right\rangle_Q \quad \forall (v,q_\odot)\in V_h \times Q_h.
\label{eq:EigenBrezzi}
\end{equation}
Then, $\lambda \geq 0$ and $\beta \approx \beta_h = \sqrt{ \lambda^{\min{}}}$.
\end{lemma}
In the framework of problem \eqref{eq:weak_P}, the eigenproblem can be recasted in the following form
involving the Schur complement of $\mathcal{A}$:
Find $p_{\odot}\in Q_h$, $\lambda_B>0$ such that\footnote{
Though it is $\lambda^2_B$ that is the eigenvalue of \eqref{eq:schur_eigw} we
shall in the following, with a slight abuse of notation, refer to
$\max \lambda_B=\lambda^{\max}_B$, $\min\lambda_B = \lambda^{\min}_B$ as the
eigenvalue bounds, extremal eigenvalues or simply eigenvalues. This choice is
intended to simplify notation and avoid proliferation of $\sqrt{}$ in the text.
Similar convention will be applied also with the Steklov eigenvalue problems
\eqref{eq:Steklov} and 
\eqref{eq:twoD_twoD_cont_mean_steklov}.
}
\begin{equation}\label{eq:schur_eigw}
BA^{-1}B'p_{\odot} = \lambda^2_B Np_{\odot}\quad\text{ in } Q'_h.
\end{equation}
Here the subscript $B$, which stands for ``Brezzi'', has been introduced
for clarity of notation as will become clear soon.
In particular, for stability of the discrete problem \eqref{eq:weak_P},
$\lambda^{\min}_{B}=\min \lambda_B$, $\lambda^{\max}_{B}=\max \lambda_B$ shall be independent
from the mesh size $h$. We remark that in this framework $\lambda^{\min}_{B}\approx\beta$.
\\ \\
\textbf{\textit{Trace constant evaluation.}} Let us consider the eigenvalue
problem defined on the domain in \Cref{fig:cont_u}, which reads as follows:
\begin{equation} \label{eq:Steklov}
\begin{aligned} 
    -\Delta u &= 0 &\quad\text{ in }\Omega\subset\mathbb{R}^2,\\
    \nabla u\cdot\nu &= \lambda_S^{-2}u &\quad\text{ on }\Gamma_{\eps},\\
    u &= 0 &\quad\text{on }\partial\Omega, 
\end{aligned} 
\end{equation}
where $\nu$ is the unit normal to $\Gamma_{\eps}$. We remark that \eqref{eq:Steklov}
is a variant of the Steklov eigenvalue problem. For rigorous mathematical
treatment of \eqref{eq:Steklov} we refer to e.g. \cite{hersch1974some} and
references therein.

The weak formulation of  \eqref{eq:Steklov}, leads to the generalized eigenproblem:
Find $u\in H^{1}_{0}(\Omega)$, $\lambda_S > 0$ satisfying
\begin{equation}\label{eq:twoD_twoD_cont_steklov}
  \int_{\Gamma_{\eps}} u v = \lambda_S^2 \int_{\Omega} \nabla u\cdot \nabla v\quad\forall v\in H^{1}_{0}
\end{equation}
so that, $\lambda^{\max}_S$ allows for estimates of the $L^2(\Gamma_{\eps})$-norm of $u$
in terms of the $H^1_0(\Omega)$-norm. Indeed, we have $\lVert u \rVert_{L^2(\Gamma_{\eps})}\leq \lambda^{\max}_S \lVert u \rVert_{H^1_0}(\Omega)$
for all $u\in {H^1_0}(\Omega)$. 

Once the eigenproblems \eqref{eq:EigenBrezzi}, \eqref{eq:Steklov} are properly discretized
(from now on subscript $h$ will denote the discretized analogues of a given quantity;
subscript $\eps_i$ will denote a quantity evaluated on a domain with inner/inclusion radius 
equal to $\eps_i$), the numerical tools for the evaluation of the inf-sup constant
$\beta$ and the trace inequality constant $C_T$ are readily established. We can, thus, proceed
with the investigation of the influence of $\eps$ in the considered mathematical framework. 

\begin{figure}
    \centering
    \includegraphics[width=1\textwidth]{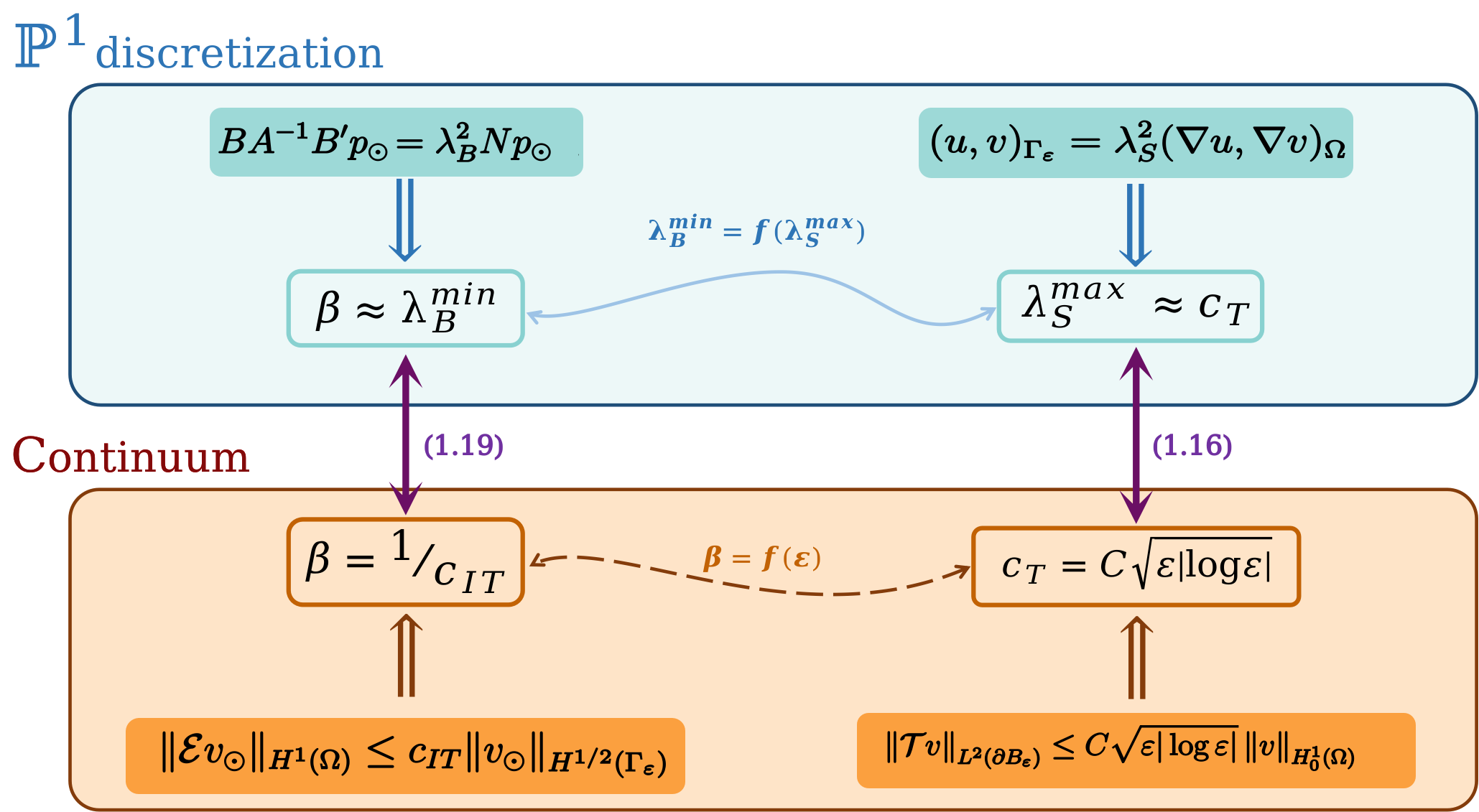}
    \vspace{-10pt}
    \caption{ Pictorial representation of the path followed in the analysis of
      the role of $\eps$ on $\beta$. Thanks to the the evaluation of the eigenvalues of problems
      \eqref{eq:schur_eigw} and \eqref{eq:Steklov}, a relation can be established between $\lambda_B^{\min{}}$ and
      $\lambda_S^{\max{}}$. Then, by virtue of equations \eqref{eq:trace_constant} \eqref{eq:twoD_twoD_cont_steklov} (represented by the purple arrows) bridging the discrete setting to the continuous one, we can close the loop and determine a relation between the Brezzi inf-sup constant $\beta$ and the inner radius $\eps$.
    }
    \label{fig:logical}
\end{figure}

The research path is summarized in \Cref{fig:logical}. Due to the analytical difficulty in the evaluation of the role of the inner radius in the inverse trace constant $C_{IT}$, and its relationship with the trace constant, we shift our inquiry to a discretized setting, by dint of the discretized eigenproblems. What will be done is a concomitant evaluation of $\{ \beta_{h,\eps_i} \}_{\eps_i \in I}$ and $\{C_{{T}_{h,\eps_i}}\}_{\eps_i \in I}$, where $I = \{\eps_0, \eps_1,\, ...,\, \eps_n \}$ and $\eps_0 > \eps_1 > ... >\eps_n $ with $\eps_n \ll 1$. Then, if some noticeable correlation between $\lambda_B^{\min}$ and $\lambda_S^{\max}$ exists, then it also holds for $\{\beta_{h,\eps_i}\}$ and $\{C_{T_{h,\eps_i}}\}$ (by means of \eqref{eq:EigenBrezzi}, \eqref{eq:schur_eigw}) so that we can find numerical evidence to support our claim. 

\subsection{Numerical results about the $2d$-$1d$ formulation}

Here we summarize the numerical experiments concerning the evaluation of $\{\beta_{h,\eps_i}\}$ and
$\{C_{T_{h,\eps_i}}\}$ with varying
$\eps_i \in I=\{10^{-1}, ..., 10^{-5}\}$. For simplicity we shall at first
  limit the investigations to $\Omega=\{x\in \mathbb{R}^2, \, |x| < 1 \}$, a circular inclusion $\Gamma_{\eps}=\left\{x\in\Omega,\, \lvert x \rvert=\eps\right\}$ and 
  $\mathbb{P}_1$ discretization.

\begin{figure}
  \centering
  \includegraphics[height=0.345\textwidth]{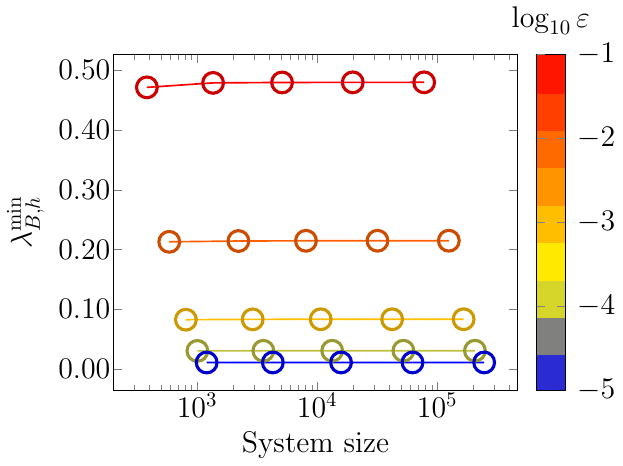}
  \hspace{17pt}
  \includegraphics[height=0.345\textwidth]{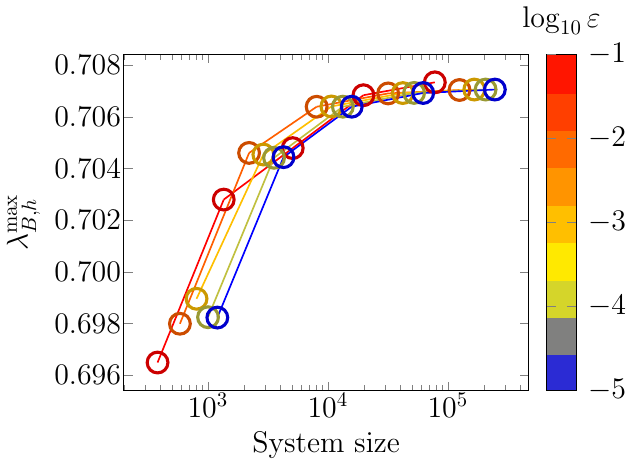}
  \vspace{-5pt}
  \caption{
    (Left) Mesh convergence of the extremal eigenvalues of Schur complement
    \eqref{eq:schur_eigw} for $\Omega=\{x\in \mathbb{R}^2,\, |x| < 1 \}$, $\Gamma_{\eps}=\left\{x\in\Omega,\, \lvert x \rvert=\eps\right\}$.
    Both spaces $V$ and $Q$ are discretized by $\mathbb{P}_1$ elements.
    For each $\eps$ a sequence of problems is considered on
    uniformly refined meshes starting from size $h_0\geq h_l\geq h_{\min}$ leading to eigenvalues
    $\lambda_{B, h_l}$. 
  }
  \label{fig:twoD_oneD}
\end{figure}

In \Cref{fig:twoD_oneD} we demonstrate that our
choice of norms and finite element spaces yields a stable problem. More
precisely, we observe that both the upper and lower bound on the Schur complement spectra are
stable in mesh refinement for each $\eps$ fixed. Moreover, $\lambda^{\max}_{B}$ appears to be bounded
in the radius of $\Gamma_{\eps}$, cf. \Cref{th:bnb}. On the other hand, $\lambda^{\min}_{B}$ seems to decrease
together with $\eps$.  

\begin{figure}
  \centering
  \includegraphics[height=0.35\textwidth]{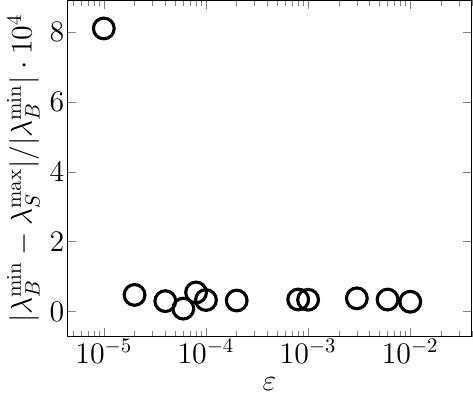}
  \hspace{17pt}
  \includegraphics[height=0.35\textwidth]{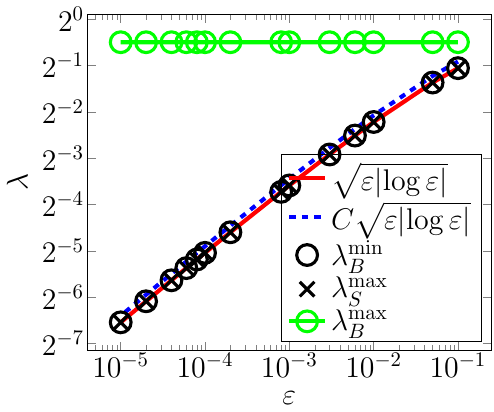}
  \vspace{-5pt}
  \caption{
    (Left)
    Error between the smallest eigenvalue $\lambda^{\min}_B$ of the Schur complement problem
    \eqref{eq:schur_eigw} and the largest eigenvalue $\lambda^{\max}_S$ of the Steklov problem
    \eqref{eq:twoD_twoD_cont_steklov}. In both cases, values from the finest level
    of refinement are considered, i.e. $\lambda_X:=\lambda_{X, h_{\min}}$.
    (Right) Dependence of the Schur complement eigenvalues on the radius of
    coupling curve $\Gamma_{\eps}=\left\{x\in\Omega, \lvert x \rvert=\eps\right\}$. Value $C\approx 0.999$ is
    obtained by fitting values $\lambda^{\min}_B$ for $\eps < 10^{-1}$.
  }
  \label{fig:twoD_oneD_cont_eps}
\end{figure}

We now propose that $\lambda^{\min}_{B}$ is closely related to the Steklov eigenvalue problem \eqref{eq:twoD_twoD_cont_steklov}.
  To investigate the relation between $\lambda^{\max}_S$ and $\lambda^{\min}_B$, \Cref{fig:twoD_oneD_cont_eps} plots the relative error
  between the two quantities\footnote{For each $\eps$ we consider a sequence of eigenvalues $\lambda_{h_l}$ computed on meshes with sizes
    $h_l$. We terminate the sequence once the relative error between subsequent eigenvalues,
    $\lvert \lambda_{h_l} - \lambda_{h_{l+1}}\rvert/\lambda_{h_l}$, is less than $10^{-4}$.
} for varying values of $\eps$. With the relative error $\sim\!\!10^{-4}$ it appears that $\lambda^{\max}_S=\lambda^{\min}_B$.
In \Cref{fig:twoD_oneD_cont_eps} we finally plot the dependence of $\lambda^{\max}_S$, $\lambda^{\min}_B$ on the radius of $\Gamma_{\eps}$. The observed dependence agrees well with the theoretical bound \cite{KVWZ}. 

\begin{figure}
  \centering
  \includegraphics[height=0.35\textwidth]{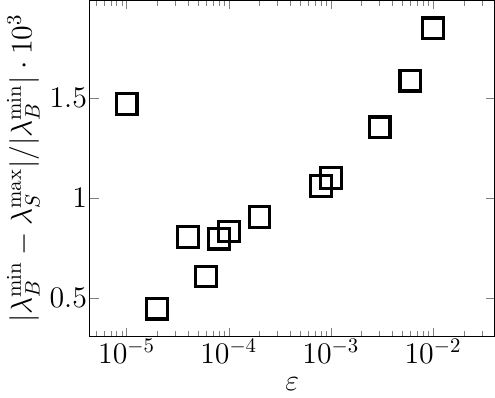}
  \hspace{17pt}
  \includegraphics[height=0.35\textwidth]{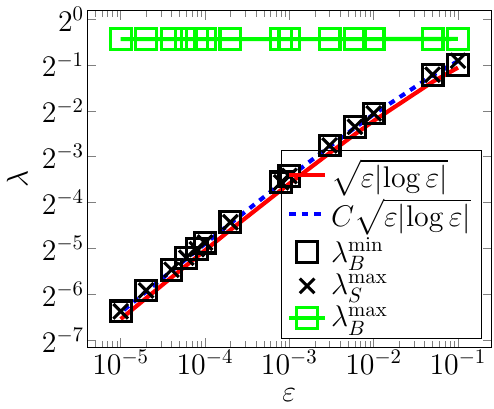}
  \vspace{-5pt}
  \caption{
    (Left)
    Error between the smallest eigenvalue $\lambda^{\min}_B$ of the Schur complement problem
    \eqref{eq:schur_eigw} and the largest eigenvalue $\lambda^{\max}_S$ of the Steklov problem
    \eqref{eq:twoD_twoD_cont_steklov} on squared domain $\Omega=(-1,1)^2$. In both cases, values from the finest level
    of refinement are considered, i.e. $\lambda_X:=\lambda_{X, h_{\min}}$.
    (Right) Dependence of the Schur complement eigenvalues on the radius of
    coupling curve $\Gamma_{\eps}=\partial(-\eps, \eps)^2$. Value $C\approx 1.116$ is
    obtained by fitting values $\lambda^{\min}_B$ for $\eps < 10^{-1}$.
  }
  \label{fig:twoD_oneD_cont_eps_square}
\end{figure}

%

In order to corroborate the independence of the established relation between $\lambda_B^{\min}$ and $\lambda_S^{\max}$
from geometrical factors, we carry out the above analysis for a square-shaped inclusion, i.e. $\Gamma_{\eps}=\partial(-\eps, \eps)^2$
in a squared domain $\Omega=(-1,1)^2$.
A comparison between the obtained values of $\lambda^{\max}_S(\eps)$ and $\lambda^{\min}_B(\eps)$, is depicted in
\Cref{fig:twoD_oneD_cont_eps_square} (see also Appendix for additional results). No significant difference is reported. The relative error between the eigenvalues remains well
below the percentage point and the expression $C\sqrt{\eps|\log \eps|}$ retraces $\lambda^{\min}_B(\eps)$ with a constant $C\approx 1.116$.

\begin{figure}
  \centering
  \includegraphics[height=0.35\textwidth]{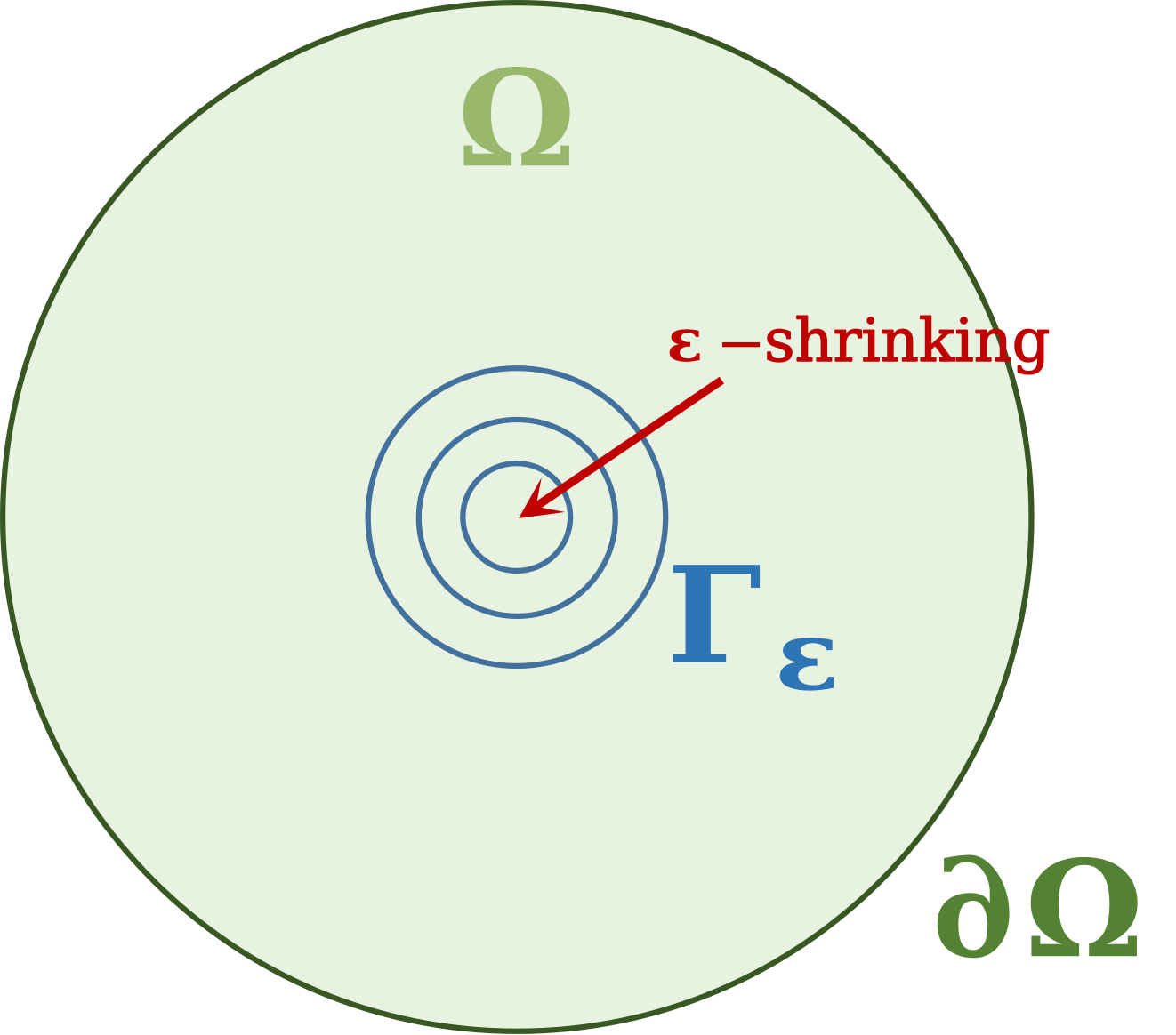}
  \hspace{27pt}
  \includegraphics[height=0.35\textwidth]{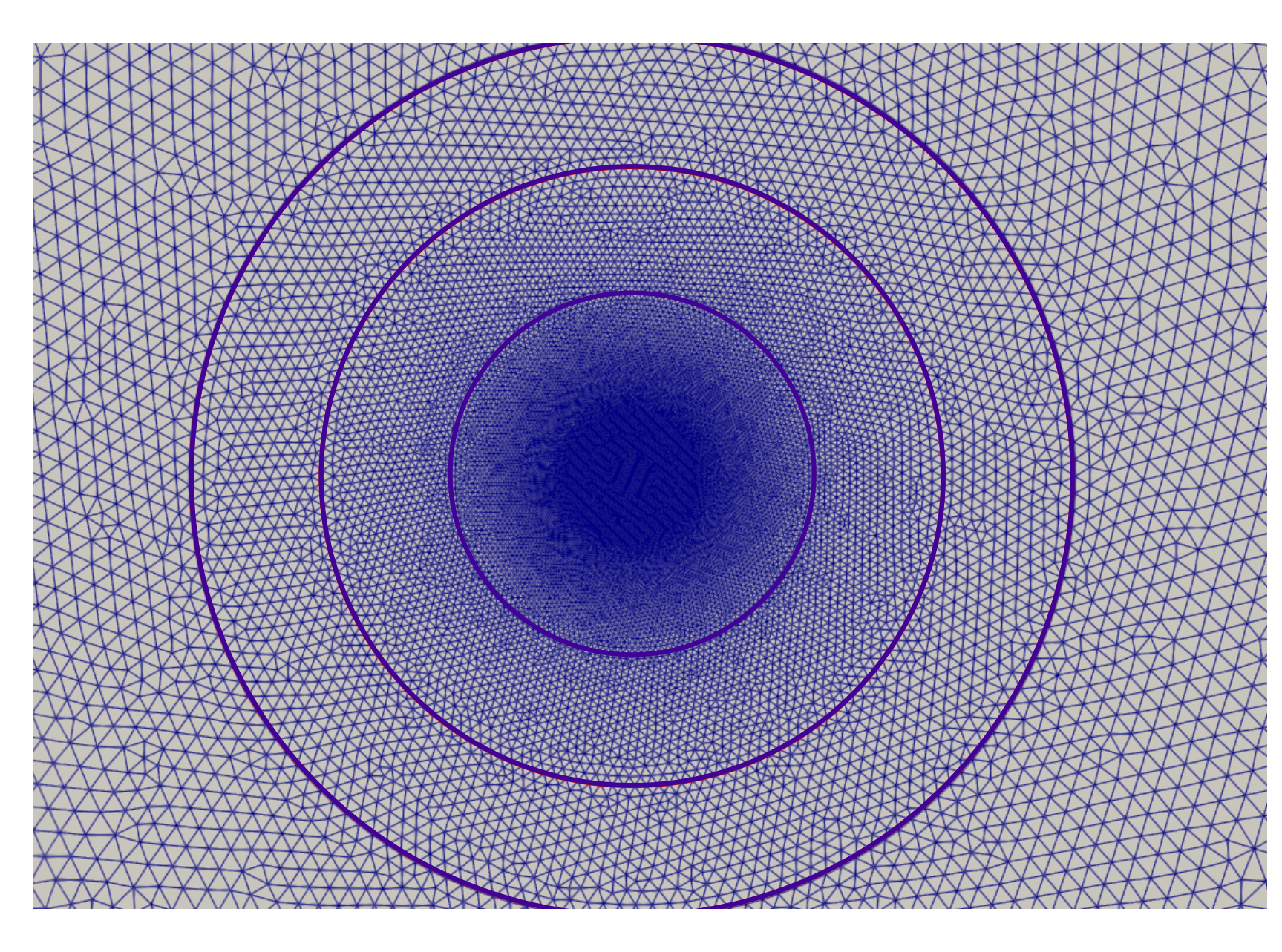}
  \vspace{-5pt}
  \caption{
    (Left) Pictorial representation of the inner-radius shrinking process. (Right) Partial
    representation (zoomed toward the $\Omega$ center) of the \emph{layered} mesh. The mesh has
    the property of being conforming simultaneously to every $\Gamma_{\eps_i}$, $\eps_i\in I$ (purple circles).
  }
  \label{fig:shrinking}
\end{figure}

\begin{figure}
  \centering
  \includegraphics[height=0.35\textwidth]{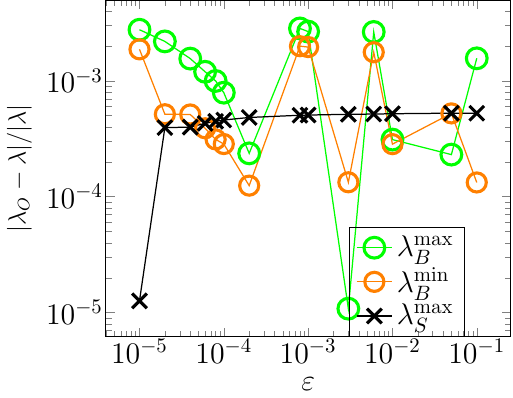}
  \hspace{10pt}
  \includegraphics[height=0.35\textwidth]{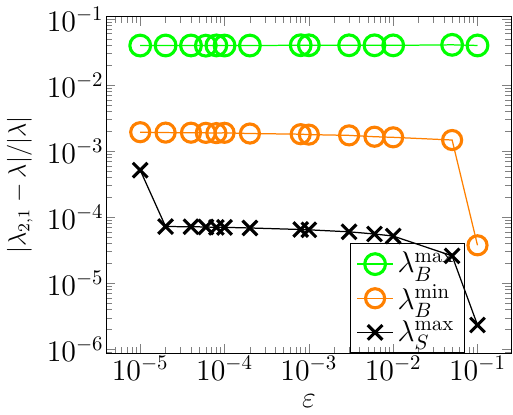}
  \vspace{-5pt}
  \caption{ 
    Relative error $|\lambda_*-\lambda|/|\lambda|$, where $\lambda$ spans $\{\lambda_{B}^{\max},\, \lambda_{B}^{\min},\, \lambda_{S}^{\max}\}$ eigenvalues obtained with standard meshing procedure (each grid conforming to a specific $\Gamma_{\eps_i}$; see \Cref{fig:twoD_oneD_cont_eps}). (Left) $\lambda_* \equiv \lambda_O$ spans $\{\lambda_{B,O}^{\max},\, \lambda_{B,O}^{\min},\, \lambda_{S,O}^{\max}\}$ obtained on the \emph{layered} mesh in \Cref{fig:shrinking} and denoted by the subscript $"O"$. (right) $\lambda_* \equiv \lambda_{2,1}$ spans $\{\lambda_{B,\,2,1}^{\max},\, \lambda_{B,\,2,1}^{\min},\, \lambda_{S,\,2,1}^{\max}\}$ obtained with a $\mathbb{P}_2-\mathbb{P}_1$ discretization and denoted by subscript $"2,1"$.
     }
  \label{fig:onion-P2P1}
\end{figure}

Having shown independence of our observations from the shape of the inclusion/coupling surface, the
effect of meshing strategy and finite element discretization will be investigated using a circular embedded domain 
$\Gamma_{\eps}$. To exclude any relevant mesh influence on the obtained results, the same analysis has
been repeated on a specific type of mesh, which we shall refer to as \emph{layered}, and which has the 
peculiarity of being $\Gamma_{\eps_{i}}$-conformal ($\Omega$ vertices and edges match with
$\Gamma_{\eps_{i}}$ ones) for every value of $\eps_i \in I$ simultaneously, cf. \Cref{fig:shrinking}. 
In such a way, the mesh configuration, which is the same for every value of the inner radius, cannot be
held responsible for any effects on $\lambda$-$\eps$ relations. The results are plotted in
\Cref{fig:onion-P2P1} (left). The relative difference between the values of $\lambda_B^{\max}(\eps)$, $\lambda_B^{\min} (\eps)$
and $\lambda_S^{\max} (\eps)$  obtained on different meshes (each conformal to a specific
$\Gamma_{\eps_i}$) and the ones obtained on the layered mesh differs by less than one percentage point. Accordingly, the same relation plotted in \Cref{fig:twoD_oneD_cont_eps} between $\lambda_S^{\max}(\eps)$ and $\lambda_B^{\min} (\eps)$
holds also for the values obtained on the layered mesh. In conclusion the results do not undergo
significant grid influence. Additional numerical results in support of this conclusion can be found in Appendix.

In addition to the geometrical (shape, mesh) factors, our aim is to exclude possible effects of
different polynomial order between the $2d$ and $1d$ spaces. The results of the comparison
between $\mathbb{P}_1-\mathbb{P}_1$  and $\mathbb{P}_2-\mathbb{P}_1$ discretization
(continuous $\mathbb{P}_2$ elements for $\Omega$ and $\mathbb{P}_1$ elements for
$\Gamma_{\eps_i}$) are plotted in \Cref{fig:onion-P2P1}. No noticeable difference can be attributed
to the change of polynomial order, except for effects on $\lambda_B^{\max}(\eps)$, which differ more
then $1\%$, but still remain independent from $\eps_i$. From this, we can deduce that also for the
$\mathbb{P}_2-\mathbb{P}_1$ discretization, $\lambda_B^{\min} \approx \lambda_S^{\max}  \sim\sqrt{\eps|\log\eps|}$.
Thus no relevant role of the polynomial degree on the results can be experienced (for detailed results see Appendix).
%

\subsection{The $2d$-$0d$ formulation for the perforated domain problem}
The link between stability of \eqref{eq:cont_u} and inner radius $\eps$
via the constant of the inverse trace inequality, which we demonstrated above
for $2d$-$1d$ trace problem, naturally extends to $3d$-$2d$ coupling as well. 
In connection with the $3d$-$1d$ problem \eqref{eq:3D1Dvar} we may then ask
if the dimensional reduction removes to observed affect of $\eps$. To address
this question we continue our investigations by applying model reduction to
\eqref{eq:cont_u} leading (formally) to the system
\begin{equation}\label{eq:twoD_twoD_cont_mean}
  \begin{aligned}
    -\Delta u &= 0&\quad\text{ in }\Omega\subset\mathbb{R}^2,\\
    \bar{u} &= g&\quad\text{ on }\Gamma,\\
    u &= 0&\quad\text{on }\partial\Omega.
  \end{aligned}
\end{equation}
We recall that $\bar{u}$ computes a mean of function $u:\Omega\rightarrow\mathbb{R}$ over
the curve $\Gamma=\Gamma_{\eps}$, i.e. $\bar{u}=\lvert\Gamma\rvert^{-1}\int_{\Gamma} u$.
Note that through $u$ and $\bar{u}$ the coupled problem \eqref{eq:twoD_twoD_cont_mean}
includes a dimensional gap of 2 (as in the case of $3d$-$1d$ problem \eqref{eq:3D1Dvar}).

\begin{figure}
  \centering
  \includegraphics[height=0.25\textwidth]{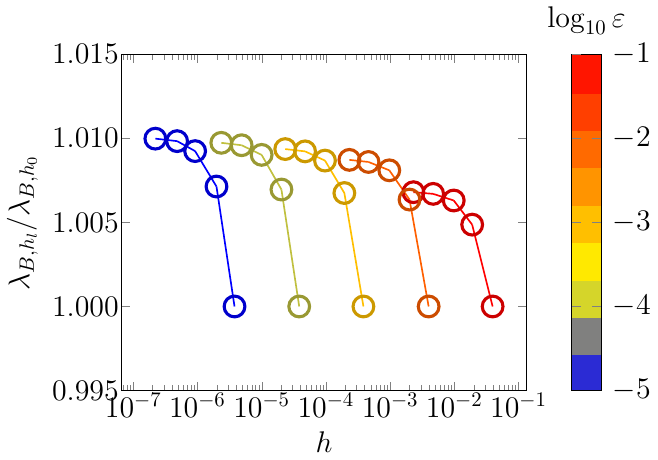}
  \includegraphics[height=0.25\textwidth]{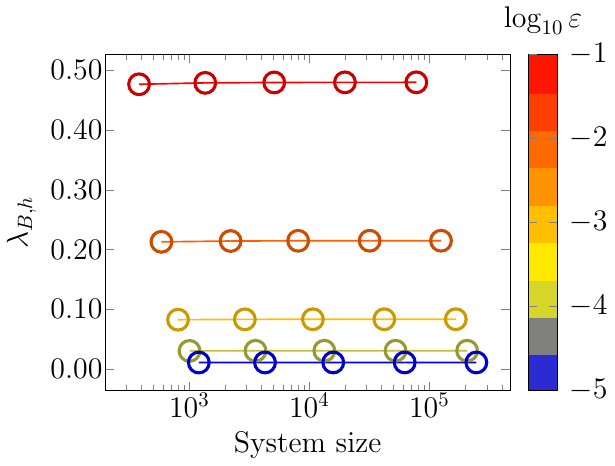}
  \includegraphics[height=0.25\textwidth]{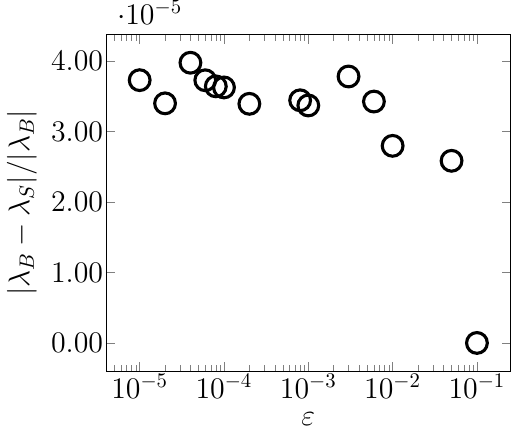}
  \vspace{-15pt}  
  \caption{
    (Left, center) Mesh convergence of the Schur complement eigenvalue of \eqref{eq:twoD_twoD_cont_mean_weak}
    for $\Omega=\{x\in \mathbb{R}^2,\, |x| < 1 \}$, $\Gamma_{\eps}=\left\{x\in\Omega, \lvert x \rvert=\epsilon\right\}$ and $\mathbb{P}_1$ elements. For
    each $\eps$ a sequence of problems is considered on uniformly refined meshes starting from size
    $h_0\geq h_l\geq h_{\min}$ leading to eigenvalues $\lambda_{B, h_l}$. (Right) Error between
    the Schur complement eigenvalues and the eigenvalues of the Steklov problem
    \eqref{eq:twoD_twoD_cont_mean_steklov}. In both cases results for $h_l=h_{\min}$ are shown.
  }
  \label{fig:twoD_oneD_mean}
\end{figure}

Letting $V=H^1_{0}{(\Omega)}$ and $Q=\mathbb{R}$ the weak form of \eqref{eq:twoD_twoD_cont_mean}
reads: Find $u\in V$, $p_\odot\in Q$ such that
\begin{equation}\label{eq:twoD_twoD_cont_mean_weak}
  \begin{aligned}
    \int_{\Omega}\nabla u\cdot \nabla v - \int_{\Gamma}p_\odot\bar{v} &= 0\quad\forall v\in V,\\
   -\int_{\Gamma}q_\odot\bar{u}                                     &= 0\quad\forall q_\odot\in Q.
  \end{aligned}
\end{equation}
For well-posedness of \eqref{eq:twoD_twoD_cont_mean_weak} when $\Gamma\subset\partial\Omega$ we
refer to \cite{Formaggia2018}. Here we shall consider $Q$ with
an inner product $(p_\odot, q_\odot)_Q=\int_{\Gamma}p_\odot q_\odot$ inducing the norm $\lVert p_\odot \rVert_Q=(p_\odot, p_\odot)^{1/2}_Q$.
On $V$, the norm is given by the $H^1$-seminorm. From the point of Brezzi theory a convenient property
of the reduced problem is the fact that $Q$ is one-dimensional. Thus, the Schur complement spectrum \eqref{eq:schur_eigw}
contains just a single eigenvalue $\lambda_B$. Then, \Cref{fig:twoD_oneD_mean} shows a computational
evidence for stability of \eqref{eq:twoD_twoD_cont_mean_weak} in $V\times Q$ with
the chosen norms. Specifically, taking $\Omega=\{x \in \mathbb{R}^2 : |x| < 1\}$ such that the triangulation of the domain
always conforms to $\Gamma_{\eps}$, and using $\mathbb{P}_1$ elements, it can be seen that
the Brezzi constant(s) related to the coupling operator\footnote{The Brezzi conditions for the bilinear
form $a$ stemming from \eqref{eq:twoD_twoD_cont_mean_weak} are easily verified and the related constants take the value 1.}
are bounded in the mesh size. However, similar to the full problem \eqref{eq:cont_u} we observe
that the numerical inf-sup/$B$-boundedness constant
$\lambda_B$ decreases with radius $\eps$, i.e. $\lambda_B=\lambda_B(\eps)$.

As with the $2d$-$1d$ problem we claim that $\lambda_B$ is closely linked to a Steklov eigenvalue.
In this case we consider: Find $u\in V$, $\lambda_S>0$ satisfying
\begin{equation}\label{eq:twoD_twoD_cont_mean_steklov}
  \int_{\Gamma}\bar{u}\bar{v} = \lambda_S^2 \int_{\Omega} \nabla u\cdot \nabla v\quad\forall v\in V.
\end{equation}
From \eqref{eq:twoD_twoD_cont_mean_steklov} it follows that the maximal eigenvalue
$\lambda_S$ relates to the estimates for the mean value of $u$ on $\Gamma_{\eps}$,
that is,
\[
\lVert \bar{u} \rVert_Q \leq \lambda_S \lVert u \rVert_V\quad\forall u\in V.
\]

The relation between the two eigenvalues is demonstrated in \Cref{fig:twoD_oneD_mean} which shows
the relative error between the eigenvalues $\lambda_B$ of the Schur complement of \eqref{eq:twoD_twoD_cont_mean_weak} and
$\lambda_S$. Here only the values obtained on the finest meshes for each $\eps$
are considered, i.e. $\lambda_S=\lambda_{S, h_{\min}}$ and analogously for $\lambda_B$.
In all the cases the observed error is $\sim\!\!10^{-3}$. 

\begin{figure}
  \centering
  \includegraphics[height=0.33\textwidth]{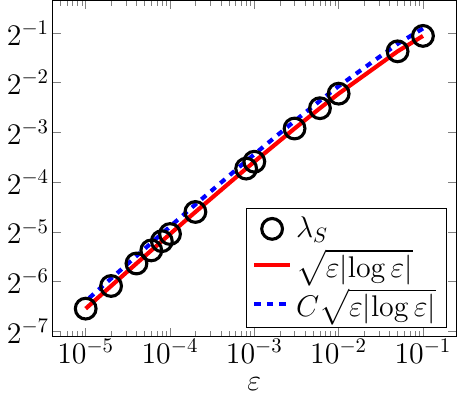}
  \hspace{15pt}
  \includegraphics[height=0.33\textwidth]{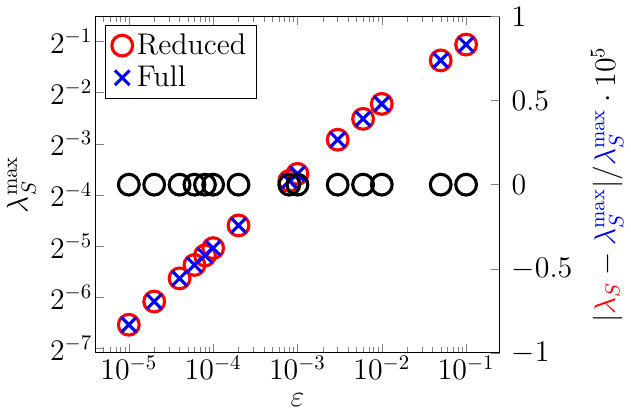}
  \vspace{-5pt}
  \caption{
    (Left) Dependence of the Steklov eigenvalue 
    $\lambda_{S, h_{\min}}$ (black $\circ$ markers) in \eqref{eq:twoD_twoD_cont_mean_steklov} on the radius of
    coupling curve $\Gamma_{\eps}=\left\{x\in\Omega, \lvert x \rvert=\epsilon\right\}$. The constant $C\approx 0.999$
    was optimized for best fit of the data.
    (Right) Comparison between the largest eigenvalues of the Steklov problems.
    The $2d$-$1d$ problem \eqref{eq:weak_P} is related to \eqref{eq:twoD_twoD_cont_steklov} with eigenvalues $\lambda^{\max}_S=\lambda^{\max}_{S, h_{\min}}$ (blue $\times$ markers), cf. \Cref{fig:twoD_oneD_cont_eps}, 
    while $\lambda_S$ (red $\circ$ markers) denotes eigenvalues of \eqref{eq:twoD_twoD_cont_mean_steklov} related to 
    $2d$-$0d$ problem \eqref{eq:twoD_twoD_cont_mean_weak}, see \Cref{fig:twoD_oneD_mean}. Relative error between the values of the full and reduced models is plotted against the right vertical axis in black $\circ$ markers.
  }
  \label{fig:twoD_oneD_mean_eps}
\end{figure}

Finally, in \Cref{fig:twoD_oneD_mean_eps} we measure the dependence of $\lambda_S$ (and $\lambda_B$)
on the radius $\eps$. It can be seen that the relation is practically identical to that
of the unreduced $2d$-$1d$ problem \eqref{eq:cont_u}, see also \cite{KVWZ}. In particular,
$\lambda_B$ goes to $0$ together with $\eps$ and  the affect of inner radius
on stability is not removed by model reduction.

In summary, we have established,
through both numerical experiments and analytical expressions about the inf-sup and trace constant, 
a clear relationship between $\beta$ and $\eps$
\begin{equation}
    \beta=C \sqrt{\eps|\log \eps|}\quad \text{as} \quad \eps \rightarrow 0
\end{equation} 
with $C$ a constant independent from $\eps$.
The result appears independent from the discretization parameters of the considered numerical framework, so that we can
assert with reasonable confidence that such behaviour is inherent in the trace/extension operator
structure of the interface problem.   

\section{Conclusion}\label{sec:conclude}

The solvability of mixed-dimensional problems plays a crucial role in effectively applying these models to real-world scenarios, such as microcirculation. In our research, we have focused on operator preconditioning as a means to address this issue. We have demonstrated that by employing suitable weighted norms, the operator preconditioning framework can successfully handle material parameters like diffusivity in both the $3d$ and $1d$ domains. However, when dealing with interface-coupled systems, these norms alone are insufficient to ensure robustness concerning geometric parameters, such as the inner radius. Through extensive numerical experiments, we have highlighted the significant impact of the parameter $\eps$ on the mathematical structure of the problem and its adverse effect on the well-posedness through the trace operator. It is worth noting that this behavior persists even in a non topologically-reduced framework. Therefore, the reduction of dimensionality and the use of appropriately scaled Sobolev spaces currently fail to guarantee the robustness of preconditioners as $\eps$ approaches zero. In our view, these findings strongly advocate for a thorough and fundamental analysis of the trace operator's role in coupling conditions within mixed-dimensional approaches. The ultimate goal is to develop a generalized trace operator capable of facilitating a robust coupling of partial differential equations across high-dimensional gaps.




\bibliographystyle{abbrv}
\bibliography{references.bib}
\clearpage
\appendix
\section{Appendix}
\subsection{Numerical experiments for square-shaped inclusion}\label{sec:appendix}

For the sake of comparison of the effect of the inclusion shape on the relation
between well-posedness of \eqref{eq:weak_P} and the diameter of the inclusion we collect
here additional results for $\Gamma_{\eps}=\partial(-\eps, \eps)^2$ and $\Omega=(-1, 1)^2$. These results are analogous
to \Cref{fig:twoD_oneD} and \Cref{fig:twoD_oneD_mean_eps} where circular $\Gamma_{\eps}$
is considered.

\begin{figure}
  \centering
  \includegraphics[height=0.345\textwidth]{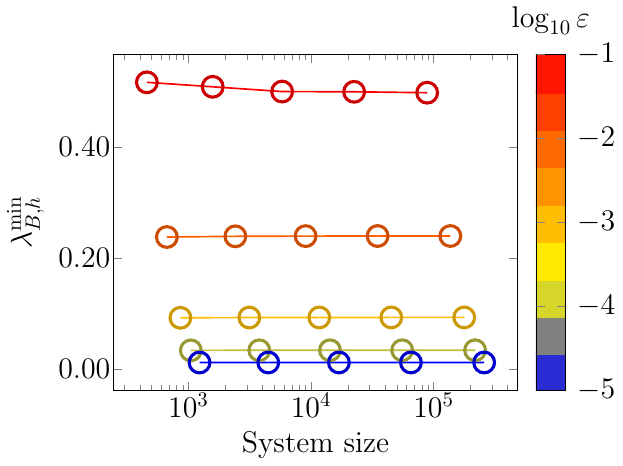}
  \hspace{15pt}
  \includegraphics[height=0.345\textwidth]{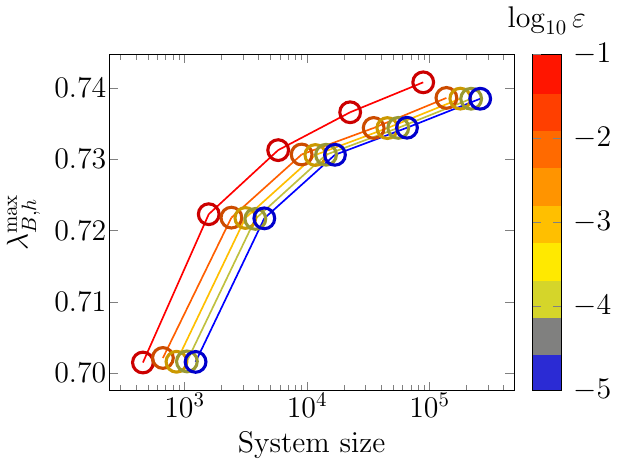}
  \vspace{-5pt}
  \caption{
    (Left) Mesh convergence of the extremal eigenvalues of Schur complement
    \eqref{eq:schur_eigw} for $\Omega=(-1, 1)^2$, $\Gamma_{\eps}=\partial(-\eps, \eps)^2$.
    Both spaces $V$ and $Q$ are discretized by $\mathbb{P}_1$ elements.
    For each $\eps$ a sequence of problems is considered on
    uniformly refined meshes starting from size $h_0\geq h_l\geq h_{\min}$ leading to eigenvalues
    $\lambda_{B, h_l}$. 
  }
  \label{fig:twoD_oneD_square}
\end{figure}

\begin{figure}
  \centering
  \includegraphics[height=0.33\textwidth]{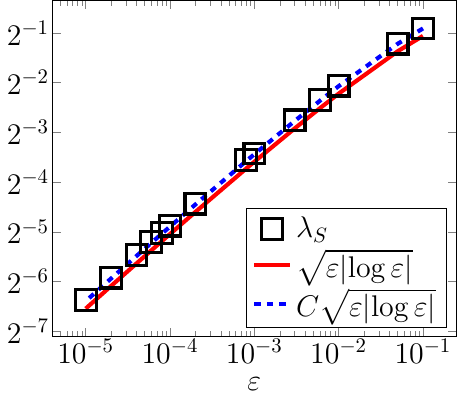}
  \hspace{15pt}
  \includegraphics[height=0.33\textwidth]{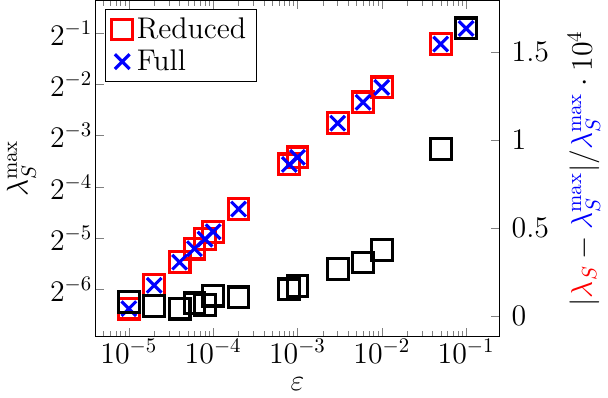}
  \vspace{-5pt}
  \caption{
    Dependence of the Steklov eigenvalue 
    $\lambda_{S, h_{\min}}$ in \eqref{eq:twoD_twoD_cont_mean_steklov} on the radius of
    coupling curve $\Gamma_{\eps}=\partial(-\eps, \eps)^2$. The constant $C\approx 1.116$
    was optimized for best fit of the data.
    (Right) Comparison between the largest eigenvalues of the Steklov problems.
    The $2d$-$1d$ problem \eqref{eq:weak_P} is related to \eqref{eq:twoD_twoD_cont_steklov} with eigenvalues $\lambda^{\max}=\lambda^{\max}_{S, h_{\min}}$, cf. \Cref{fig:twoD_oneD_cont_eps}, 
    while $\lambda_S$ denotes eigenvalues of \eqref{eq:twoD_twoD_cont_mean_steklov} related to 
    $2d$-$0d$ problem \eqref{eq:twoD_twoD_cont_mean_weak}, see \Cref{fig:twoD_oneD_mean}.
Relative error between the values of the full and reduced models is plotted against the right vertical axis in black $\square$ markers.  }
  \label{fig:twoD_oneD_mean_eps_square}
\end{figure}

 \subsection{Numerical experiments for layered mesh}
Here we collect additional results regarding the numerical experiments for a circular domain $\Omega=\{x\in \mathbb{R}^2,\, |x| < 1 \}$ with circular inclusion $\Gamma_{\eps}=\left\{x\in\Omega,\, \lvert x \rvert=\eps\right\}$ obtained on layered mesh (see \Cref{fig:shrinking} (right)) conformal simultaneously to every $\Gamma_{\eps_i}$, with $\eps_i \in \{10^{-1}, 10^{-2}, 10^{-3}, 10^{-4}, 10^{-5} \} $.
\begin{figure}[h]
  \centering
  \includegraphics[height=0.35\textwidth]{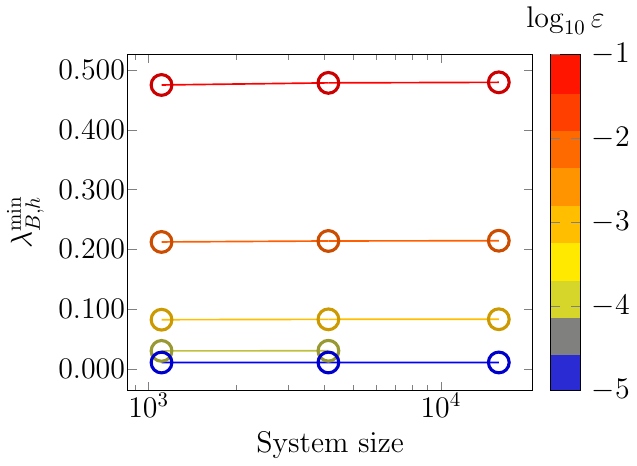}
  \hspace{10pt}    
  \includegraphics[height=0.35\textwidth]{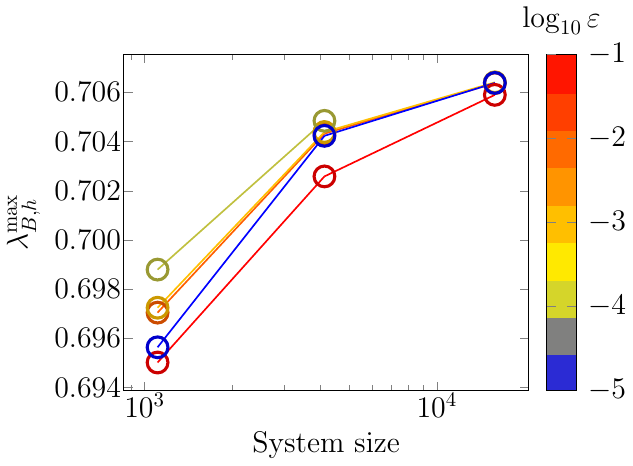}
  \vspace{-5pt}
  \caption{
    (Left) Mesh convergence of the extremal eigenvalues of Schur complement
    \eqref{eq:schur_eigw} for $\Omega=\{x\in \mathbb{R}^2,\, |x| < 1 \}$, $\Gamma_{\eps}=\left\{x\in\Omega,\, \lvert x \rvert=\eps\right\}$.
    Both spaces $V$ and $Q$ are discretized by $\mathbb{P}_1$ elements.
    For each $\eps$ a sequence of problems is considered on
    uniformly refined layered mesh starting from size $h_0\geq h_l\geq h_{\min}$ leading to eigenvalues
    $\lambda_{B, h_l}$.
    }
  \label{ }
\end{figure}

\begin{figure}
  \centering
  \includegraphics[height=0.35\textwidth]{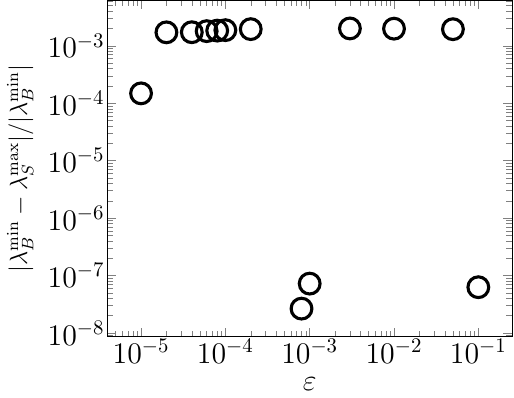}
  \hspace{15pt}    
  \includegraphics[height=0.36\textwidth]{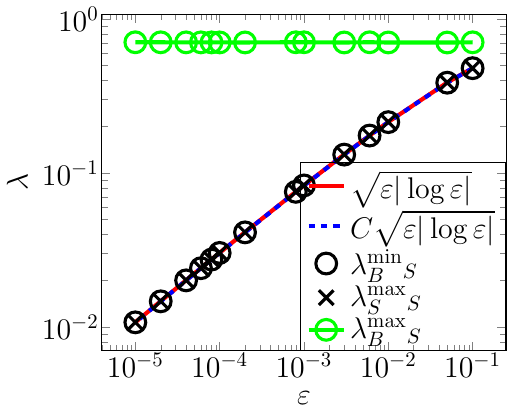}
  \vspace{-15pt}
  \caption{
    Error between the smallest eigenvalue $\lambda^{\min}_B$ of the Schur complement problem
    \eqref{eq:schur_eigw} and the largest eigenvalue $\lambda^{\max}_S$ of the Steklov problem
    \eqref{eq:twoD_twoD_cont_steklov}. In both cases, values from the finest level
    of refinement are considered, i.e. $\lambda_X:=\lambda_{X, h_{\min}}$.
    (Right) Dependence of eigenvalues from the radius of
    coupling curve $\Gamma_{\eps}=\left\{x\in\Omega, \lvert x \rvert=\eps\right\}$. Value $C\approx 0.999$ is
    obtained by fitting values $\lambda^{\min}_B$ for $\eps < 10^{-1}$.
    }
  \label{fig:onion_p1p1}
\end{figure}

\subsection{Numerical experiments for $\mathbb{P}_2-\mathbb{P}_1$ discretization}
\begin{figure}
  \centering
  \includegraphics[height=0.35\textwidth]{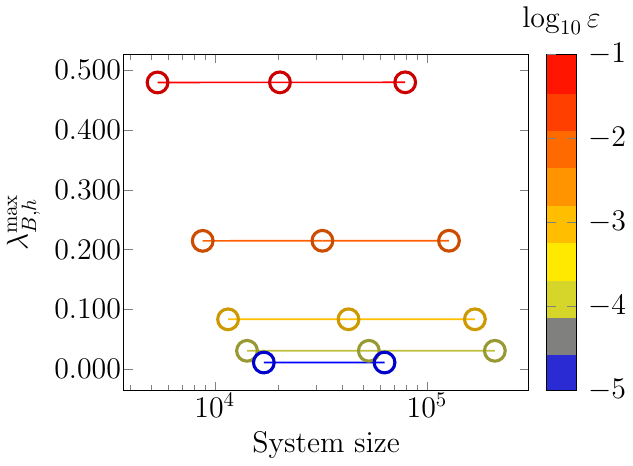}
  \hspace{10pt}  
  \includegraphics[height=0.35\textwidth]{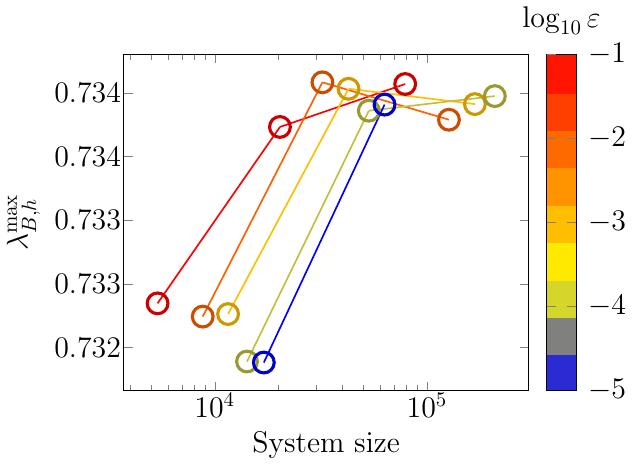}
  \vspace{-5pt}
  \caption{
    (Left) Mesh convergence of the extremal eigenvalues of Schur complement \eqref{eq:schur_eigw} for $\Omega=\{x\in \mathbb{R}^2,\, |x| < 1 \}$, $\Gamma_{\eps}=\left\{x\in\Omega,\, \lvert x \rvert=\eps\right\}$. Spaces $V$ and $Q$ are discretized respectively by $\mathbb{P}_2$ and $\mathbb{P}_1$ elements. For each $\eps$ a sequence of problems is considered on uniformly refined meshes starting from size $h_0\geq h_l\geq h_{\min}$ leading to eigenvalues
    $\lambda_{B, h_l}$.
    }
  \label{fig:space}
\end{figure}
\begin{figure}
  \centering
  \includegraphics[height=0.35\textwidth]{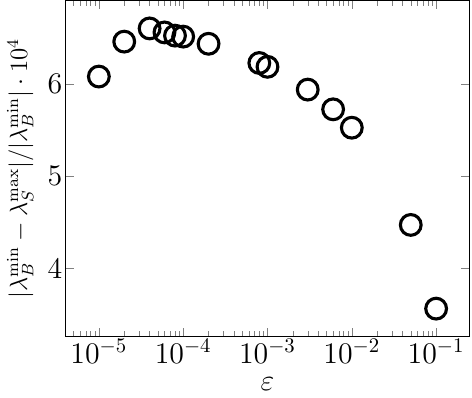}
  \hspace{15pt}
  \includegraphics[height=0.35\textwidth]{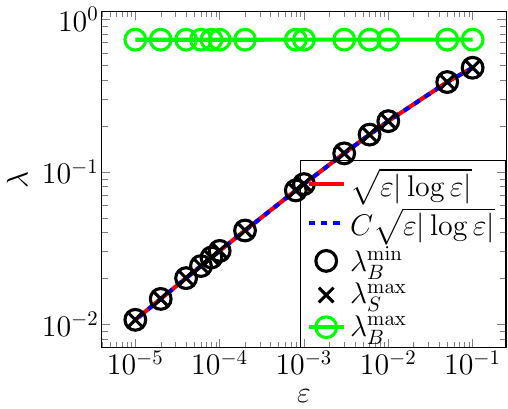}
  \vspace{-5pt}
  \caption{
    Error between the smallest eigenvalue $\lambda^{\min}_B$ of the Schur complement problem
    \eqref{eq:schur_eigw} and the largest eigenvalue $\lambda^{\max}_S$ of the Steklov problem
    \eqref{eq:twoD_twoD_cont_steklov}. In both cases, values from the finest level
    of refinement are considered, i.e. $\lambda_X:=\lambda_{X, h_{\min}}$.
    (Right) Dependence of eigenvalues from the radius of
    coupling curve $\Gamma_{\eps}=\left\{x\in\Omega, \lvert x \rvert=\eps\right\}$. Value $C\approx 0.999$ is
    obtained by fitting values $\lambda^{\min}_B$ for $\eps < 10^{-1}$. 
    }
  \label{fig:onion_p2p1}
\end{figure}
\end{document}